\documentclass[11pt]{article}

\usepackage{arxiv}
\usepackage{mathtools,amssymb,subfigure,xcolor}
\usepackage{siunitx}
\usepackage[utf8]{inputenc}
\usepackage{helvet}

\usepackage{setspace}
\usepackage[hidelinks]{hyperref}

\usepackage[utf8]{inputenc} 
\usepackage[T1]{fontenc}    
\usepackage{hyperref}       
\usepackage{url}            
\usepackage{booktabs}       
\usepackage{amsfonts}       
\usepackage{nicefrac}       
\usepackage{microtype}      
\usepackage{lipsum}
\usepackage[square,numbers]{natbib}

\title{Hamiltonian pitchfork bifurcation in transition across index-1 saddles}
\author{Wenyang Lyu\thanks{wl16298@bristol.ac.uk}, \, Shibabrat Naik\thanks{s.naik@bristol.ac.uk}, \, Stephen Wiggins\thanks{s.wiggins@bristol.ac.uk}\\
School of Mathematics, University of Bristol\\
Fry building, Woodland Road, Bristol BS8 1UG, UK
}

\date{}
\begin{document}
\maketitle


\begin{abstract}
	We study the effect of changes in the parameters of a two dimensional potential energy surface on the phase space structures relevant for chemical reaction dynamics. The changes in the potential energy are representative of chemical reactions such as isomerization between two structural conformations or dissociation of a molecule with an intermediate. We present a two degrees of freedom quartic Hamiltonian that shows pitchfork bifurcation when the parameters are varied and we derive the bifurcation criteria relating the parameters. Next, we describe the phase space structures \textemdash~unstable periodic orbits and its associated invariant manifolds, and phase space dividing surfaces \textemdash~for the systems that can show trajectories undergo reaction defined as crossing of a potential energy barrier. Finally, we quantify the reaction dynamics for these systems by obtaining the directional flux and gap time distribution to illustrate the dependence on total energy and coupling strength between the two degrees of freedom.
\end{abstract}

%


\section{Introduction}

Understanding, predicting, and controlling transition dynamics has significance in many physical, chemical, biological, and engineering systems. For example, ionization of a hydrogen atom under an electromagnetic field in atomic physics~\cite{JaFaUz2000}, transport of defects in solid state and semiconductor physics~\cite{Eckhardt1995}, chemical reaction rates and pathways in 
chemical physics~\cite{Komatsuzaki1999,WiWiJaUz2001}, dynamic buckling or collapse in structural mechanics~\cite{Collins2012,ZhViRo2018,sieber_nonlinear_2019}, ship capsize~\cite{Virgin1989,ThDe1996,naik2017geometry}, escape and recapture of comets and asteroids in celestial mechanics~\cite{JaRoLoMaFaUz2002,DeJuLoMaPaPrRoTh2005,Ross2003}, and escape into inflation or re-collapse to singularity in 
cosmology~\cite{DeOliveira2002}. In these systems, understanding transitions between stable configurations constitutes finding the phase space structures and analysing their stability. These phase space structures are equilibrium points, periodic orbits, and invariant manifolds (surfaces or curves in the phase space on which an initial condition once started stays forever) in the phase space of the governing equations of motion. In a typical scenario of transition dynamics, these phase space structures form the skeleton for the transition by guiding the trajectories between stable regions in the phase space~\cite{WiWiJaUz2001,UzJaPaYaWi2002}. This indicates that studying the changes in the stability of these phase space structures, that is bifurcation analysis, gives an understanding of the dependence of the transitions on the parameters. 

In the context of reaction dynamics, transition between stable regions in the phase space represents reaction by breaking and forming of chemical bonds. The dynamical systems perspective of chemical reactions entails describing the breaking and forming of chemical bonds using the geometry of phase space structures. The starting point for analysing the geometry and tracking any changes in the stability of these structures is the bifurcation of equilibrium points of the vector field. In the setting of Hamiltonian models where the total energy (or the Hamiltonian) is a sum of kinetic plus potential energy, the equilibrium points of the vector field are also the critical points on the potential energy surface (PES). The potential energy in the Hamiltonian is obtained from the electronic structure calculations performed under the Born-Oppenheimer approximation which separates the low frequency atomic nuclei motion from the high frequency electronic motion~\cite{wales_energy_2004}. The topographical features of a given potential energy function such as the mountain top, valley floor, valley ridge, and saddles influence the reaction in a dynamical sense, that is their presence is felt for non-zero momentum (or kinetic energy) even though they are configuration space (coordinates that define the potential energy) features. Thus, it is natural to expect that the transition dynamics is affected by the changes in the features of the PES when the parameters of the Hamiltonian model are varied; for example solvent mass in solution-phase reaction~\cite{garcia-meseguer_influence_2019}, energy of the excitation and total energy of the molecule~\cite{farantos_energy_2009}. 

The issue of bifurcation in Hamiltonian models of chemical reactions has received some attention in the context of specific reactions such as collinear $HgI_2 \rightarrow HgI + I$ reaction~\cite{burghardt_molecular_1995}, collinear and non-collinear $H_2 + H$ exchange reaction~\cite{li_bifurcation_2009,Inarrea11}, bimolecular reactions~\cite{mackay2014bifurcations}. In these studies, the phase space dividing surface~\cite{Wigner38,Wigner39,wiggins_geometry_1990} undergoes bifurcation as the total energy is increased, aside from Ref.~\cite{li_bifurcation_2009} who also considered the role of bath modes on the location of the phase space dividing surface. Ref.~\cite{mackay2014bifurcations} studied the changes in the topology of the transition state (which in general is called the NHIM for N-DOF systems) as the total energy in increased. From the dynamical systems point of view and as also present in Wigner's formulation for two DOF systems in the phase space~\cite{Wigner38,Wigner39}, the correct transition state is the unstable periodic orbit (UPO) which depends on the total energy of the system, and thus its stability can undergo bifurcation. Furthermore, as the total energy is increased the region around the saddle equilibrium point on the energy surface can deviate from the bottleneck geometry. This can manifest in loss of the \emph{locally non-recrossing}~\cite{waalkens2004direct} property required by the phase space dividing surface and thus leads to breakdown of the transition state theory~\cite{burghardt_molecular_1995,Inarrea11}. So studying the effect of total energy on reaction dynamics using completely integrable (where the degrees-of-freedom are separable or uncoupled) Hamiltonian systems explains the changes in the geometry of the phase space structures relevant to the reactions~\cite{mauguiere2013bifurcations}. Here, we will consider a two DOF Hamiltonian as a simple model of a reaction in a bath (liquid solvent) where the bifurcation of the unstable periodic orbit, called Hamiltonian pitchfork bifurcation, occurs due to changes in the parameters. 

The description of a reacting molecular system is highly influenced by the topography of the potential energy surface. Potential wells describe ``products'' and saddles are the ``transition points'' between wells. This is a general mathematical structure of the transition dynamics, and therefore it would seem reasonable to transfer the language of ``chemical transformation'' to other fields with a similar potential energy surface description. One such physical system and which is of engineering interest is the buckling of beams and columns under compression. In modeling and analysing these structural elements, transition dynamics arises in studying the motion between stable configurations that represent mode shapes of these structural elements~\cite{lawrence_phonon_2002,carr_buckling_2005,chakraborty_buckled_2011,Collins2012,sieber_nonlinear_2019,ZhViRo2018}. For example, the analysis of a dynamical model of a nanobeam under compression in Ref.~\cite{Collins2012} considers the classical  Euler-Bernoulli beam equation subject to compressive stress applied at both ends. The authors considered a two-mode truncation of the forced beam equation. More precisely, the authors considered parameter regimes where the first mode is unstable and the second mode can be either stable or unstable, and the remaining (neglected) modes are always stable. Material parameters were used corresponding to a silicon nanobeam. The two-mode model Hamiltonian is the sum of a (diagonal) kinetic energy term and a potential energy term.  The form of the potential energy function that the authors obtained was analogous with isomerization reactions in chemistry, where ``isomerization'' in this case corresponded to a transition between two stable beam configurations. Thus, the dynamics of the buckled beam was studied using the conceptual framework for the theory of isomerization reactions. When the second mode is stable the potential energy surface has an index one saddle, and when the second mode is unstable the potential energy surface has an index two saddle and two index one saddles. Symmetry of the system allowed the authors to readily construct a phase space dividing surface between the two ``isomers'' (buckled states); and they rigorously proved that, in a specific energy range, it is a normally hyperbolic invariant manifold~\cite{wiggins_normally_2014}. The energy range is sufficiently wide that we were able to  treat the effects of the index one and index two saddles on the isomerization dynamics in a unified fashion. We computed reactive fluxes, mean gap times, and reactant phase space volumes for three stress values at several different energies. In all cases the phase space volume swept out by isomerizing trajectories is considerably less than the reactant density of states, proving that the dynamics is highly non-ergodic. This was supported by the associated gap time distributions with one or more ``pulses'' of trajectories. Computation of the reactive flux correlation function showed no sign of a plateau region; rather, the flux exhibited oscillatory decay, indicating that, for the two-mode model Hamiltonian in the physical regime considered, a rate constant for ``isomerization'' (buckling) does not exist. These insights are possible by bringing the developments in methods for studying chemical reactions to dynamical buckling of structures. In this article, we analyse the Hamiltonian pitchfork bifurcation relevant for transition dynamics where changes in the parameters of the potential energy shows dynamics on a closed double-well, a closed single-well, and an open single-well system. 


In Sect.~\ref{sect:model}, we describe the Hamiltonian model and analyse the linear stability of the equilibrium points to find the bifurcation criteria relating the parameters. In Sect.~\ref{sect:results}, we discuss the phase space structures governing the transition dynamics and obtain quantitative measures of the transition using the phase space structures. Finally, in Sect.~\ref{sect:conclusions}, we present our conclusions of this work and outlook for related future work.

\section{Model Hamiltonian for the pitchfork bifurcation\label{sect:model}} 

We consider a two degrees of freedom (DOF) Hamiltonian where a quartic potential is coupled to a harmonic oscillator with a bilinear coupling. 
\begin{align}
	\mathcal{H}(x,y,p_x,p_y) = T(p_x,p_y) + V(x,y) = \dfrac{1}{2} \left( \frac{p_x^2}{m_s} + \frac{p_y^2}{m_b} \right) - \frac{\alpha  }{2}\, x^2 + \frac{\beta }{4} \, x^4 + \dfrac{\omega}{2} y^2 + \dfrac{\varepsilon}{2} \left(x-y\right)^2 
	\label{eqn:ham_2dof_bilinear_coupling}
\end{align}
This form of a model Hamiltonian is borrowed from the solution-phase reactions literature where a finite sum of harmonic oscillators are referred to as the ``bath'' coupled with the quartic Hamiltonian referred to as the ``system''. Furthermore, we define the forward transition (or forward reaction) as a trajectory crossing from $x > 0$ to $x < 0$, and the backward transition (or backward reaction) is defined as crossing from $x < 0$ to $x > 0$. 

The Hamilton's equations are given by
\begin{equation}
	\begin{aligned}
		\dot{x} &= \dfrac{\partial \mathcal{H}}{\partial p_x} =  \frac{p_x}{m_s} \\
		\dot{y} &= \dfrac{\partial \mathcal{H}}{\partial p_y} = \frac{p_y}{m_b} \\
		\dot{p_x} &= -\dfrac{\partial \mathcal{H}}{\partial x} =  \alpha \, x- \beta \, x^3 + \varepsilon (y - x)\\
		\dot{p_y} &= -\dfrac{\partial \mathcal{H}}{\partial y} = -\omega y + \varepsilon (x-y)
	\end{aligned}\label{hameq_2dof}
\end{equation}
The equilibria for this system are located at $\mathbf{x}_1^e = (0,0,0,0)$ and $\pm \mathbf{x}_2^e = \pm \left(x^e,y^e,0,0\right)$.
where
\begin{equation}
	x^e = \sqrt{ \frac{1}{\beta} \left(\alpha- \frac{\omega\varepsilon}{\omega + \varepsilon} \right)} >0 \;,\quad y^e = x^e \left( \frac{\varepsilon}{\omega +\varepsilon}  \right) >0.
	\label{eqn:cSpace_eqCoords}
\end{equation}
where we require either $\beta < 0$ and $\alpha < \omega \varepsilon / (\omega + \varepsilon)$ or $\beta > 0$ and $\alpha > \omega \varepsilon / (\omega + \varepsilon)$ for the system to have three equilibria. In the case when $\beta < 0$ and $\alpha > \omega \varepsilon / (\omega + \varepsilon)$ or $\beta > 0$ and $\alpha < \omega \varepsilon / (\omega + \varepsilon)$, the system admits one equilibrium point at $\mathbf{x}_1^e = (0,0,0,0)$.

\subsection{Linear stability and bifurcation analysis}

The Jacobian of the Hamiltonian vector field is given by

\begin{align}
	\mathbb{J}(x,y,p_x,p_y) =  \begin{pmatrix}
		\dfrac{\partial^2 \mathcal{H}}{\partial x\partial p_x} & \dfrac{\partial^2 \mathcal{H}}{\partial y\partial p_x}&  \dfrac{\partial^2 \mathcal{H}}{\partial p_x^2} & \dfrac{\partial^2 \mathcal{H}}{\partial p_y\partial p_x}\\
		\dfrac{\partial^2 \mathcal{H}}{\partial x\partial p_y} & \dfrac{\partial^2 \mathcal{H}}{\partial y\partial p_y}& \dfrac{\partial^2 \mathcal{H}}{\partial p_x\partial p_y} & \dfrac{\partial^2 \mathcal{H}}{\partial p_y^2}\\
		-\dfrac{\partial^2 \mathcal{H}}{\partial x^2} &  -\dfrac{\partial^2 \mathcal{H}}{\partial y\partial x}&
		-\dfrac{\partial^2 \mathcal{H}}{\partial p_x\partial x} & -\dfrac{\partial^2 \mathcal{H}}{\partial p_y\partial x}\\
		-\dfrac{\partial^2 \mathcal{H}}{\partial x\partial y} & -\dfrac{\partial^2 \mathcal{H}}{\partial y^2}&  -\dfrac{\partial^2 \mathcal{H}}{\partial p_x\partial y} & -\dfrac{\partial^2 \mathcal{H}}{\partial p_y\partial y}
	\end{pmatrix}
	=  \begin{pmatrix}
		0  &  0 &  1 & 0      \\
		0  &  0 &  0 & 1      \\
		\alpha -\varepsilon- 3\beta x^2  &  \varepsilon & 0 &0\\
		\varepsilon  &-\omega-\varepsilon &  0&0
	\end{pmatrix}
\end{align}

The eigenvalues of $\mathbb{J}(x,y,p_x,p_y)$ evaluated at $\mathbf{x}_1^e$ are $\pm \sqrt{\lambda_1},\pm i \sqrt{-\lambda_2}$, where 
\begin{equation}
	\lambda_{1,2} = \dfrac{\alpha - \omega- 2\varepsilon \pm \sqrt{(\omega+\alpha)^2+4\varepsilon^2}}{2}.
\end{equation}
The term inside the square root is positive which means that $\lambda_2$ is negative for $\alpha \leq 0$. $\lambda_2$ is also negative for $\alpha > 0$ as the term inside the square root is larger than or equal to $\alpha$ when $\alpha > 0$. The critical value of $\alpha$ such that $\lambda_1=0$ occurs when 
\begin{equation}
	\alpha = \frac{\omega \varepsilon}{\omega+\varepsilon}.
\end{equation}

The eigenvalues of $\mathbb{J}(x,y,p_x,p_y)$ evaluated at $\pm \mathbf{x}_2^e$ are $\pm \sqrt{\lambda_3},\pm i \sqrt{-\lambda_4}$, where 
\begin{align}
	\lambda_{3,4} 
	&= \dfrac{-2\alpha - \omega- 2\varepsilon +3\omega \varepsilon/(\omega+\varepsilon)\pm \sqrt{[2\alpha-(\omega+3\omega\varepsilon/(\omega+\varepsilon))]^2+4\varepsilon^2}}{2}\\
	&= \dfrac{-2\alpha - (\omega^2+2\varepsilon^2)/(\omega+\varepsilon)\pm \sqrt{[2\alpha-(\omega+3\omega\varepsilon/(\omega+\varepsilon))]^2+4\varepsilon^2}}{2}.
\end{align}
The term inside the square root is positive which means that $\lambda_4$ is negative for $\alpha \geq 0$. $\lambda_4$ is also negative for $\alpha < 0$ as the term inside the square root is larger than or equal to $-2\alpha$ when $\alpha < 0$. The critical value of $\alpha$ such that $\lambda_3=0$ occurs when 
\begin{equation}
	\alpha = \frac{\omega \varepsilon}{\omega+\varepsilon}.
\end{equation}
We note that the existence of the critical value of $\alpha$ depends on the sign of $\alpha$ and will be discussed in the following sections. We also get the expected result that the eigenvalues are independent of $\beta$ since the linearization of the vector field at the equilibria only depends on the quadratic terms in the Hamiltonian.

The total energy of the equilibrium points are
\begin{align}
	\mathcal{H}(\mathbf{x}_1^e) = & 0, \\
	\mathcal{H}(\pm \mathbf{x}_2^e) = & (x^e)^2 \left(-\frac{\alpha}{2} + \frac{\beta}{4}(x^e)^2+ \frac{\omega\varepsilon}{2(\omega + \varepsilon)} \right) = -\dfrac{1}{4\beta} \left(\alpha- \frac{\omega\varepsilon}{\omega + \varepsilon} \right)^2.
	\label{energy_eqpts}
\end{align}
We note that the depth of the potential well is simply the difference between the saddle equilibrium point and the center equilibrium point. In the system considered here, the depth is the energy $\dfrac{1}{4\beta} \left(\alpha - \frac{\omega\varepsilon}{\omega + \varepsilon} \right)^2$. We also denote the excess energy $\Delta E$ as the energy difference between the total energy of the system, $e$ and the energy of the saddle equilibrium point. Now, we classify the linear stability type of the equilibria.

\textbf{Case $\alpha< \varepsilon \omega/(\varepsilon+\omega)$:}
\begin{enumerate}
	\item $\beta  < 0$. Three eqilibrium points, one at $\mathbf{x}_1^e = (0,0,0,0)$. $\lambda_1<0,\lambda_2<0$ which means that $\mathbf{x}_1^e$ a centre x centre equilibrium point.
	The other two equilibrium points are at $\pm \mathbf{x}_2^e = \pm (x^e,y^e,0,0)$. $\lambda_3>0,\lambda_4<0$ which means that $\pm \mathbf{x}_2^e$ are saddle-centre equilibrium points.
	\item $\beta  \geq 0$. One equilibrium point at $\mathbf{x}_1^e = (0,0,0,0)$. $\lambda_1<0,\lambda_2<0$ which means that $\mathbf{x}_1^e$ is a centre-centre equilibrium point.
\end{enumerate}

\textbf{Case $\alpha= \varepsilon \omega/(\varepsilon+\omega)$:}
\begin{enumerate}
	\item $\beta  = 0$. Line of equilibrium points at
	\begin{equation*}
		y = x (\frac{\varepsilon}{\omega+\varepsilon})
	\end{equation*}
	\item $\beta  \neq 0$. One equilibrium point at $\mathbf{x}_1^e = (0,0,0,0)$. $\lambda_1=0,\lambda_2<0$ which means that $\mathbf{x}_1^e$ is a nonhyperbolic equilibrium point.
\end{enumerate}

\begin{figure}[!ht] 
	\centering 
	\subfigure[]{\includegraphics[width=0.32\textwidth]{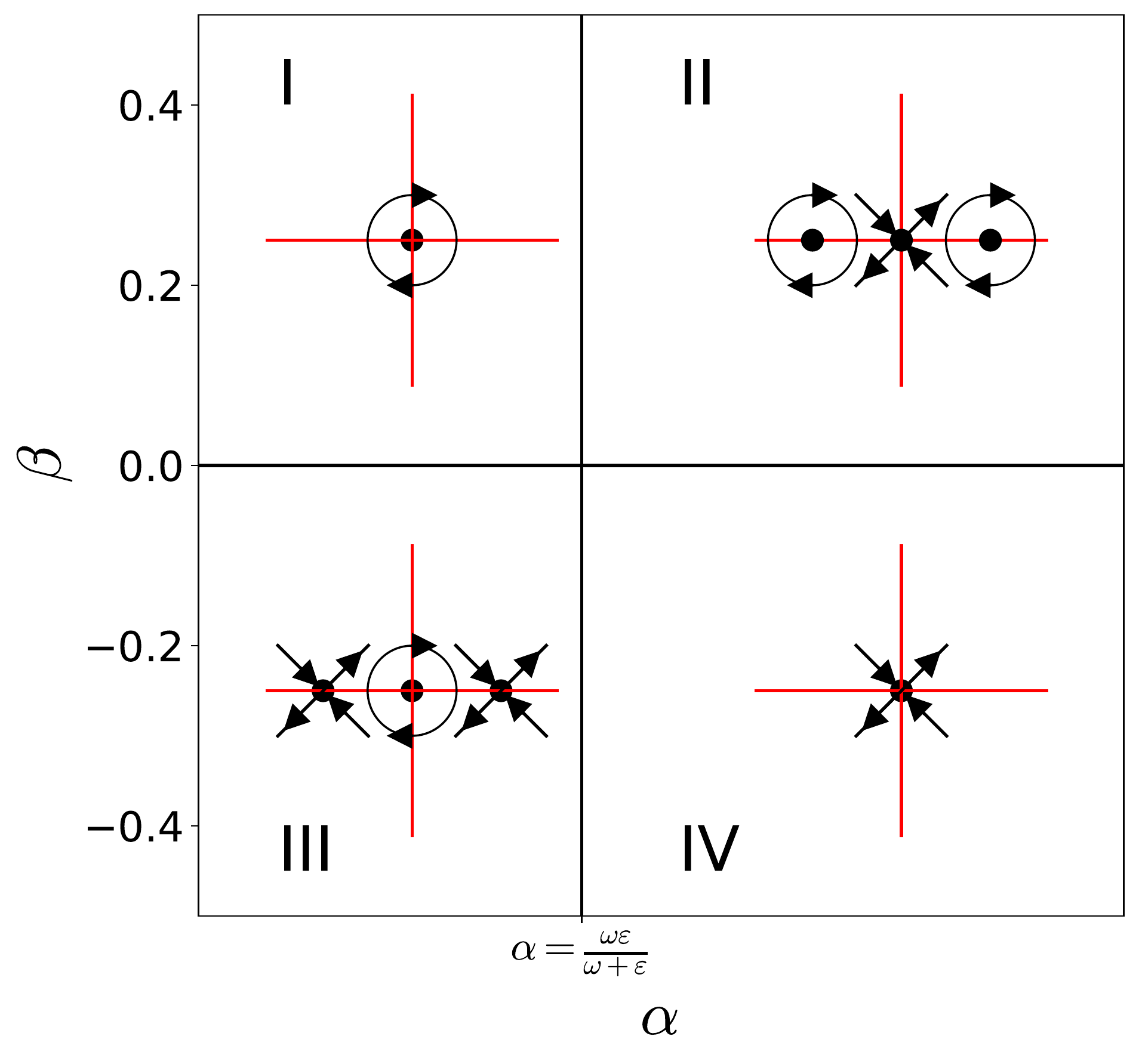}}
	\subfigure[]{\includegraphics[width=0.32\textwidth]{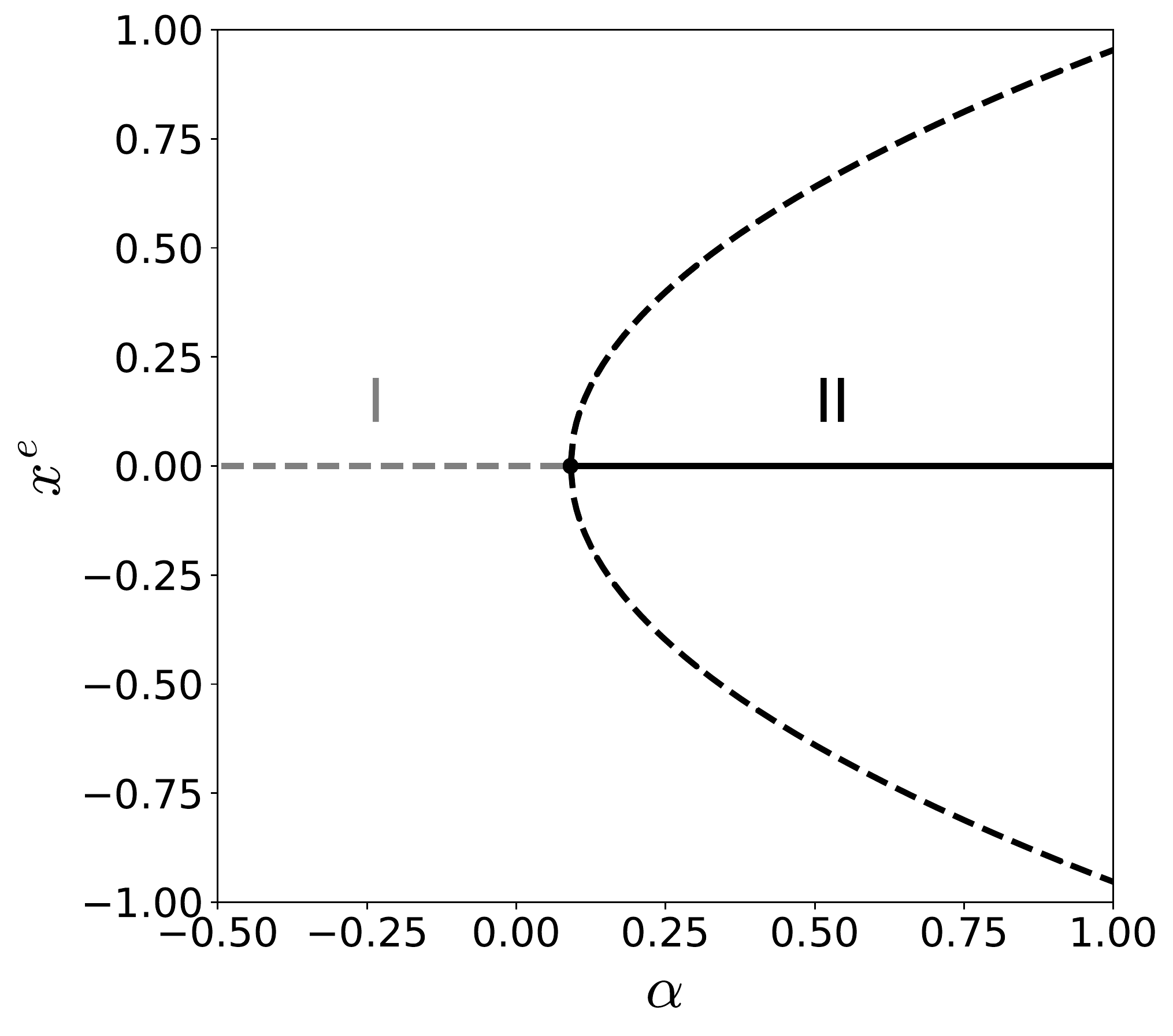}}
	\subfigure[]{\includegraphics[width=0.32\textwidth]{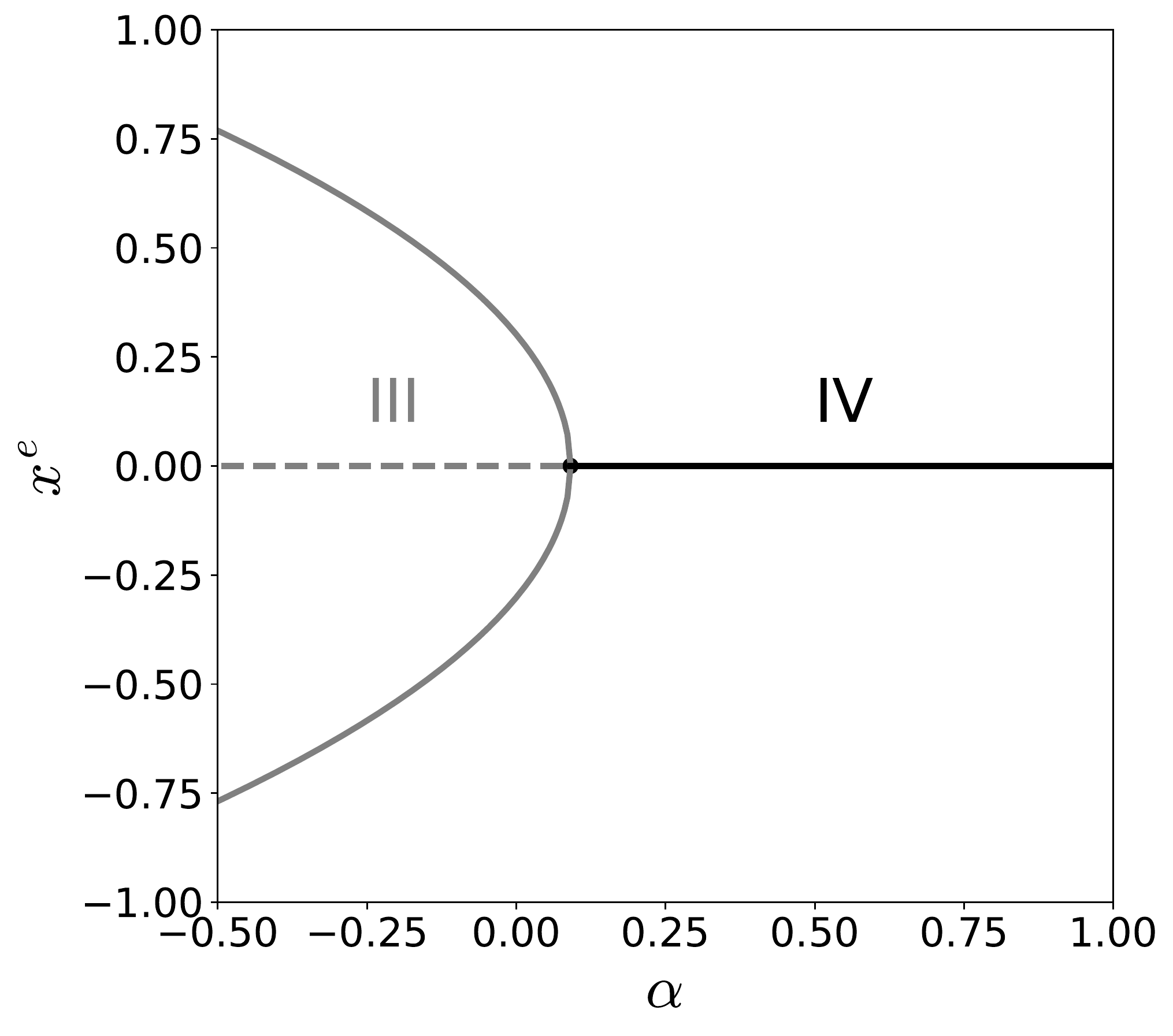}}
	\caption{(a) Bifurcation diagram for the coupled system with $\omega=1,\varepsilon=0.1$. Each subfigure correspond to the phase space geometry for 4 representative choices(I: $\alpha<\varepsilon \omega/(\varepsilon+\omega), \beta>0$, II: $\alpha>\varepsilon \omega/(\varepsilon+\omega), \beta>0$, III: $\alpha<\varepsilon \omega/(\varepsilon+\omega), \beta<0$, IV: $\alpha>\varepsilon \omega/(\varepsilon+\omega), \beta<0$) of $\alpha, \beta$ in the $x-p_x$ plane. Bifurcation diagram for the coupled system with (b) $\beta=1$, (c) $\beta=-1$. Other parameters $\omega=1,\varepsilon=0.1$.}
	\label{fig:bifurcation_diag}
\end{figure}

\textbf{Case $\alpha> \varepsilon \omega/(\varepsilon+\omega)$:}
\begin{enumerate}
	\item $\beta \leq 0$. One equilibrium point at $\mathbf{x}_1^e = (0,0,0,0)$. $\lambda_1>0,\lambda_2<0$ which means that $\mathbf{x}_1^e$ is a saddle-centre equilibrium point.
	\item $\beta  > 0$. Three eqilibrium points, one at $\mathbf{x}_1^e = (0,0,0,0)$. $\lambda_1>0,\lambda_2<0$ which means that $\mathbf{x}_1^e$ is a saddle-centre equilibrium point.
	The other two equilibrium points are at $\pm \mathbf{x}_2^e = \pm (x^e,y^e,0,0)$, with $\lambda_3<0$, $\lambda_4<0$ which means that $\pm \mathbf{x}_2^e$ are centre-centre equilibrium points.
\end{enumerate}
The summary of the bifurcation analysis of the equilibrium points associated to $\lambda_1$ or $\lambda_3$ for the coupled system is shown in Fig. ~\ref{fig:bifurcation_diag}. The linearised stability of the equilibrium points associated to $\lambda_2$ or $\lambda_4$ is not shown as it is guaranteed to have centre type stability.

\begin{figure}[!h] 
	\centering 
	\subfigure[]{\includegraphics[width=0.24\textwidth]{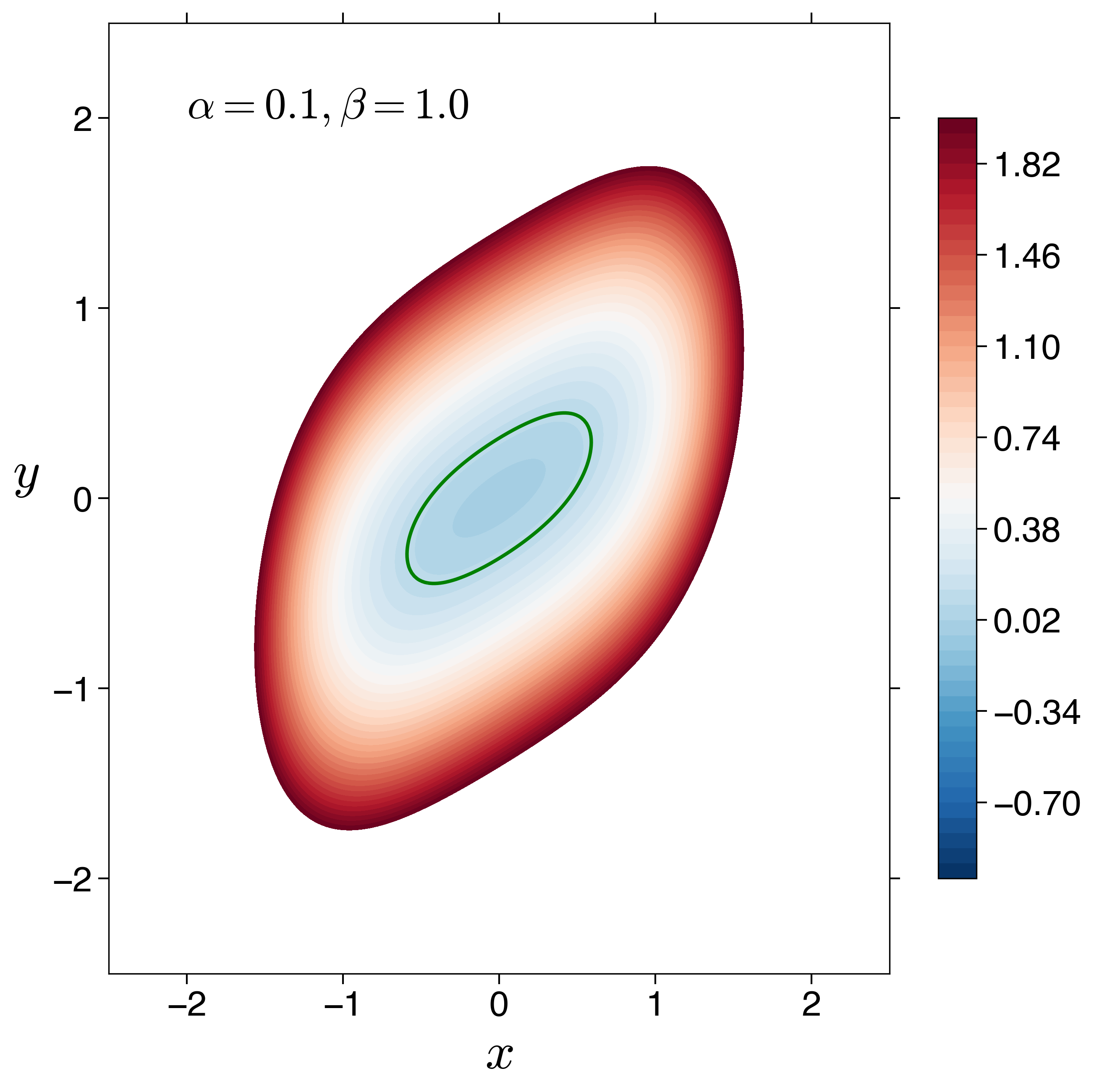}}
	\subfigure[]{\includegraphics[width=0.24\textwidth]{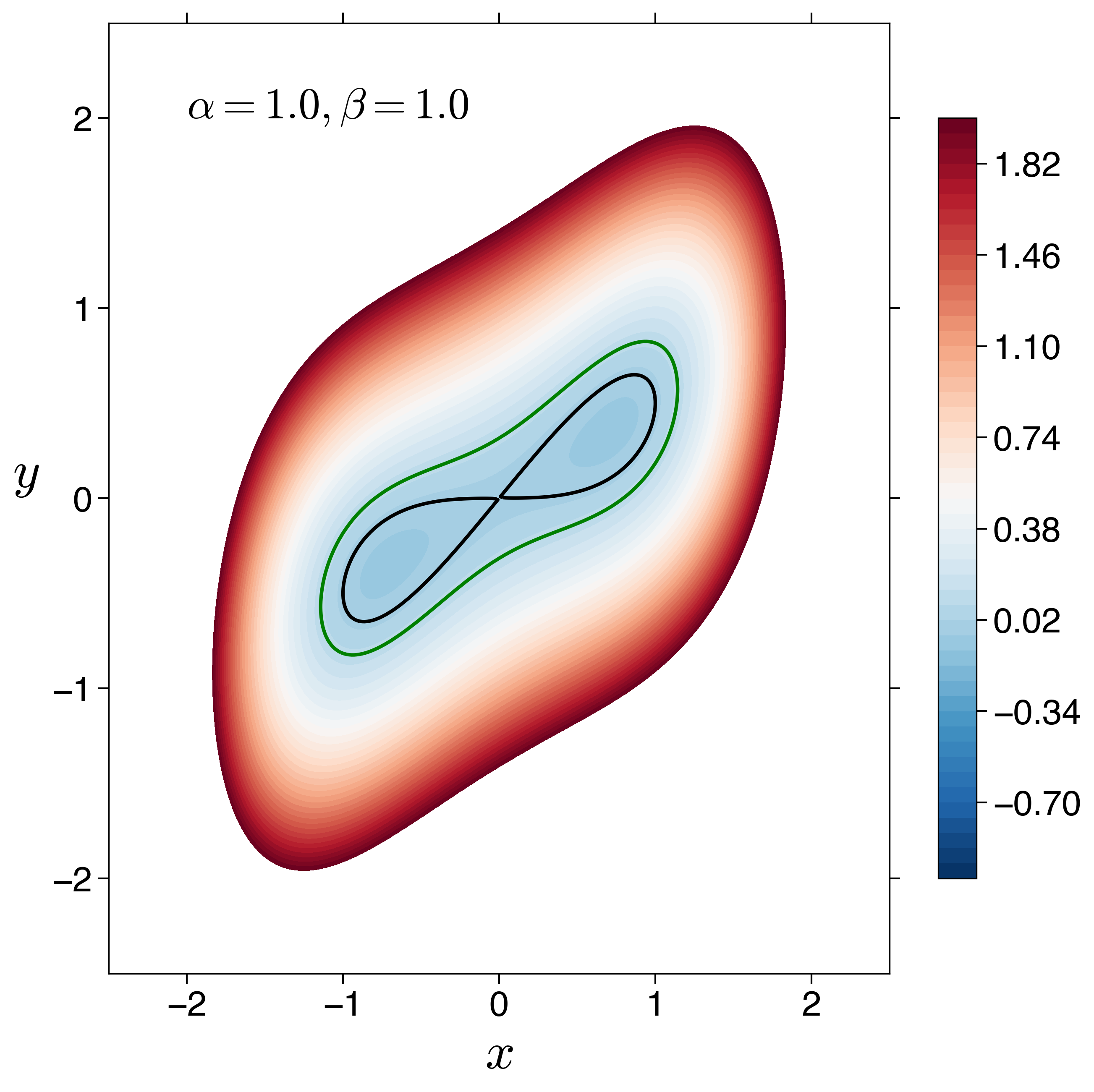}}
	\subfigure[]{\includegraphics[width=0.24\textwidth]{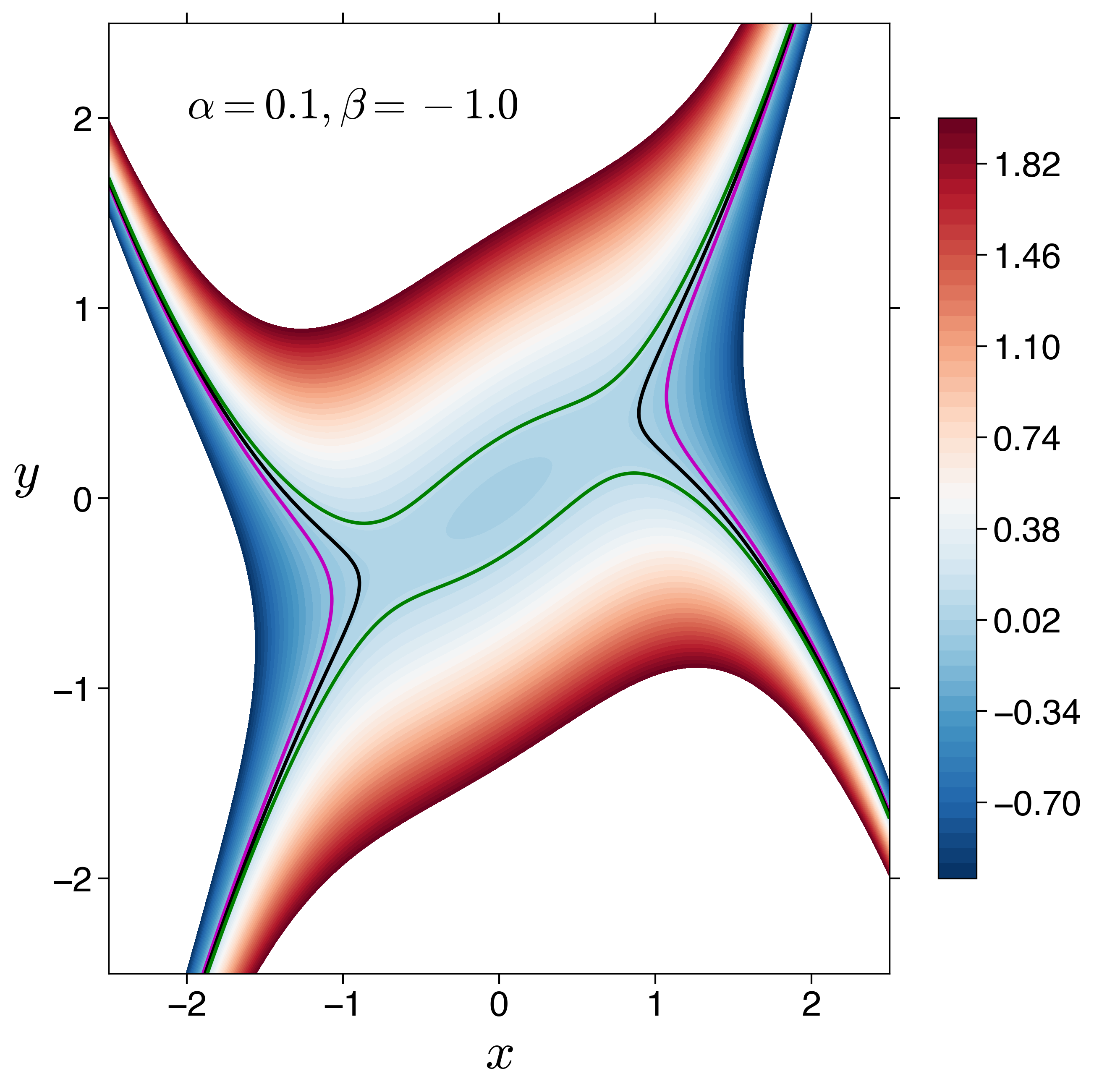}}
	\subfigure[]{\includegraphics[width=0.24\textwidth]{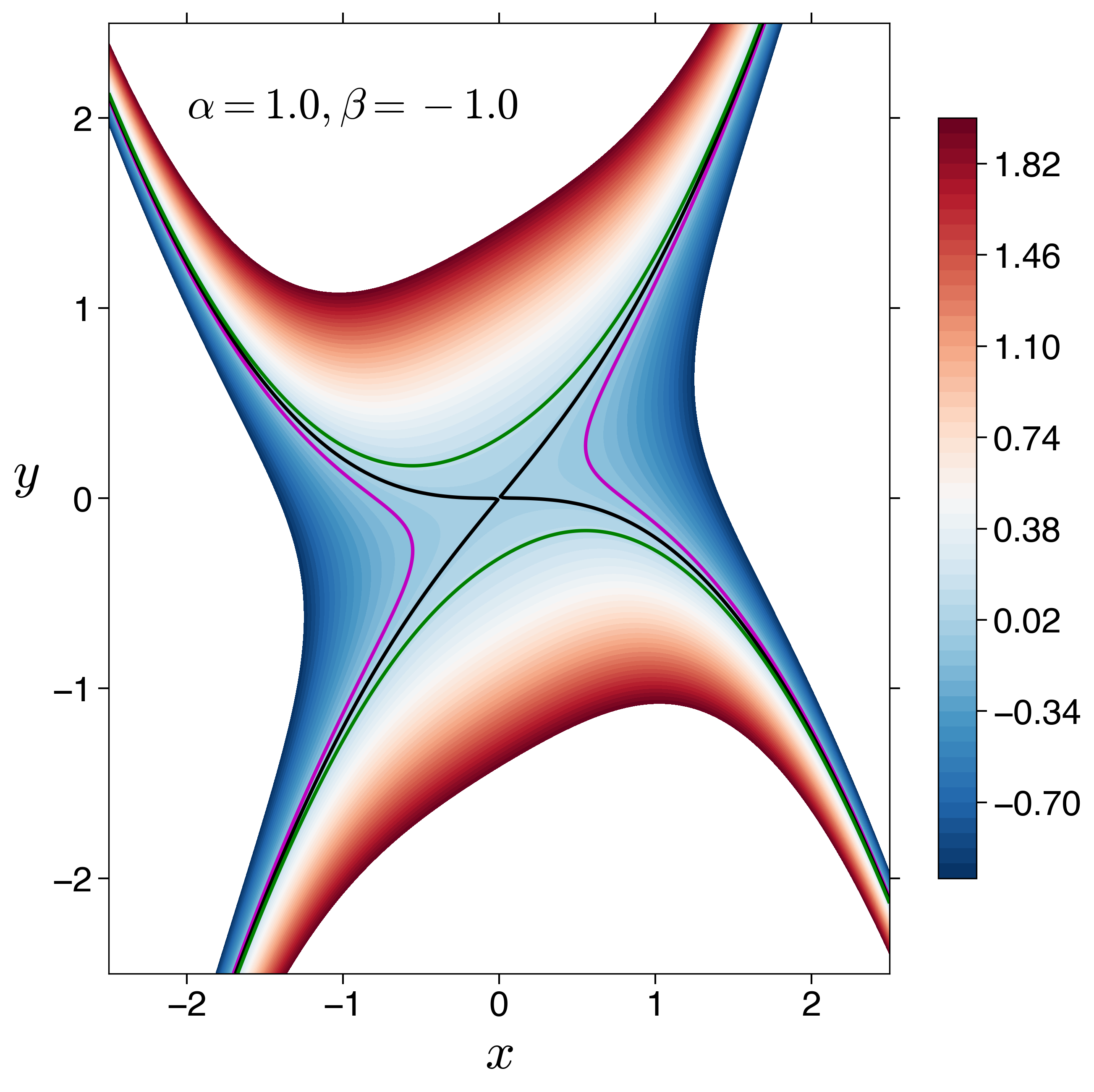}}
	\caption{Changes in the topology of the potential energy surface with changes in the parameters, $\alpha, \beta$. The green, black, and magenta curves represent the equipotential contours at total energy $0.1, 0, -0.1$, respectively. (a-d) Typical potential energy surfaces for \textbf{Cases I-IV}, respectively.}
	\label{fig:pes_contourf_bilinear}
\end{figure} 

\section{Results and discussion\label{sect:results}}

\subsection{Phase space structures relevant for the transition} 
\begin{figure}[!h] 
	\centering 
	\subfigure[$(1, 2, 0.5, 0.1)$]{\includegraphics[width=0.48\textwidth]{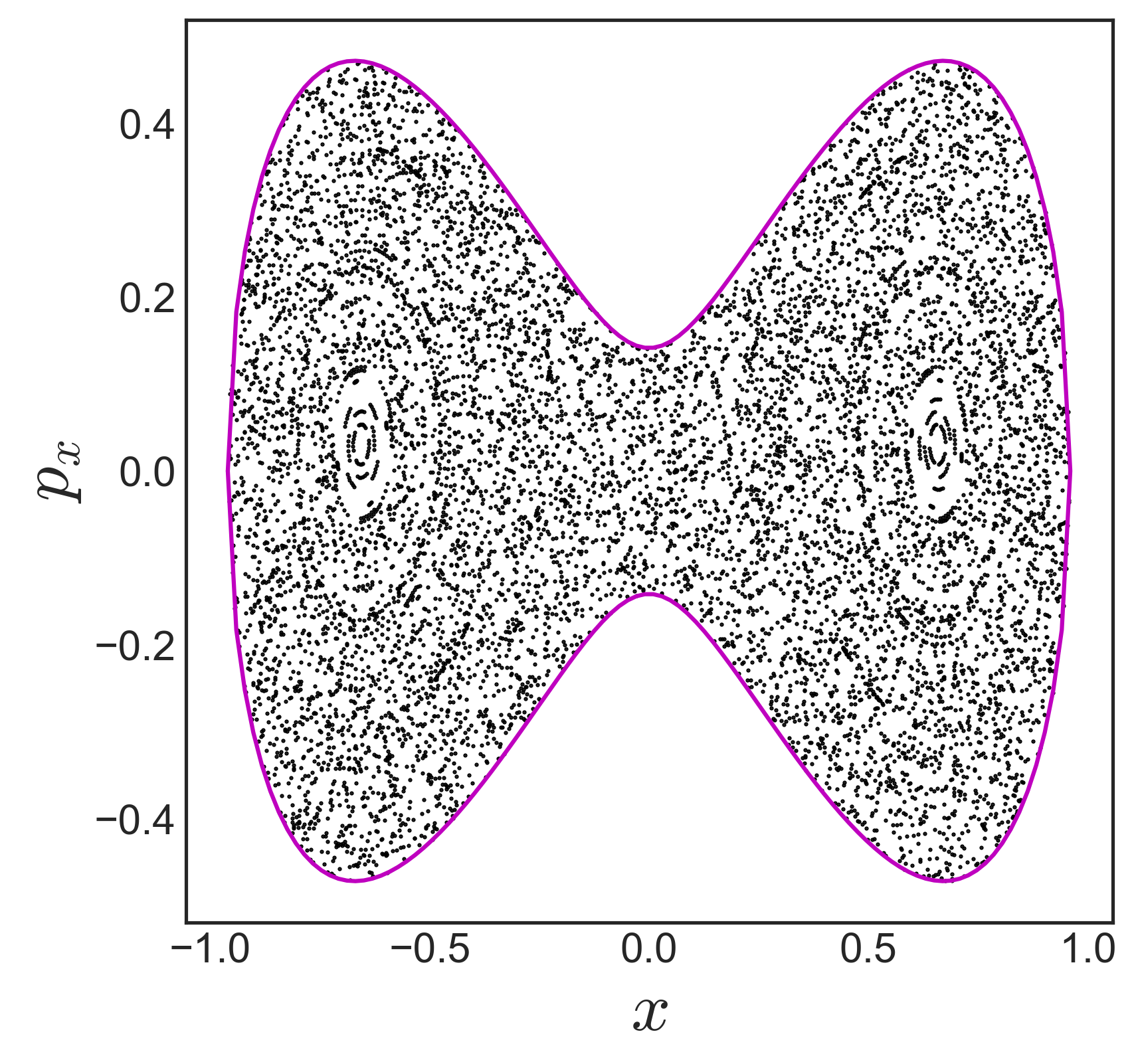}}
	\subfigure[$(1, 2, 4, 0.5)$]{\includegraphics[width=0.48\textwidth]{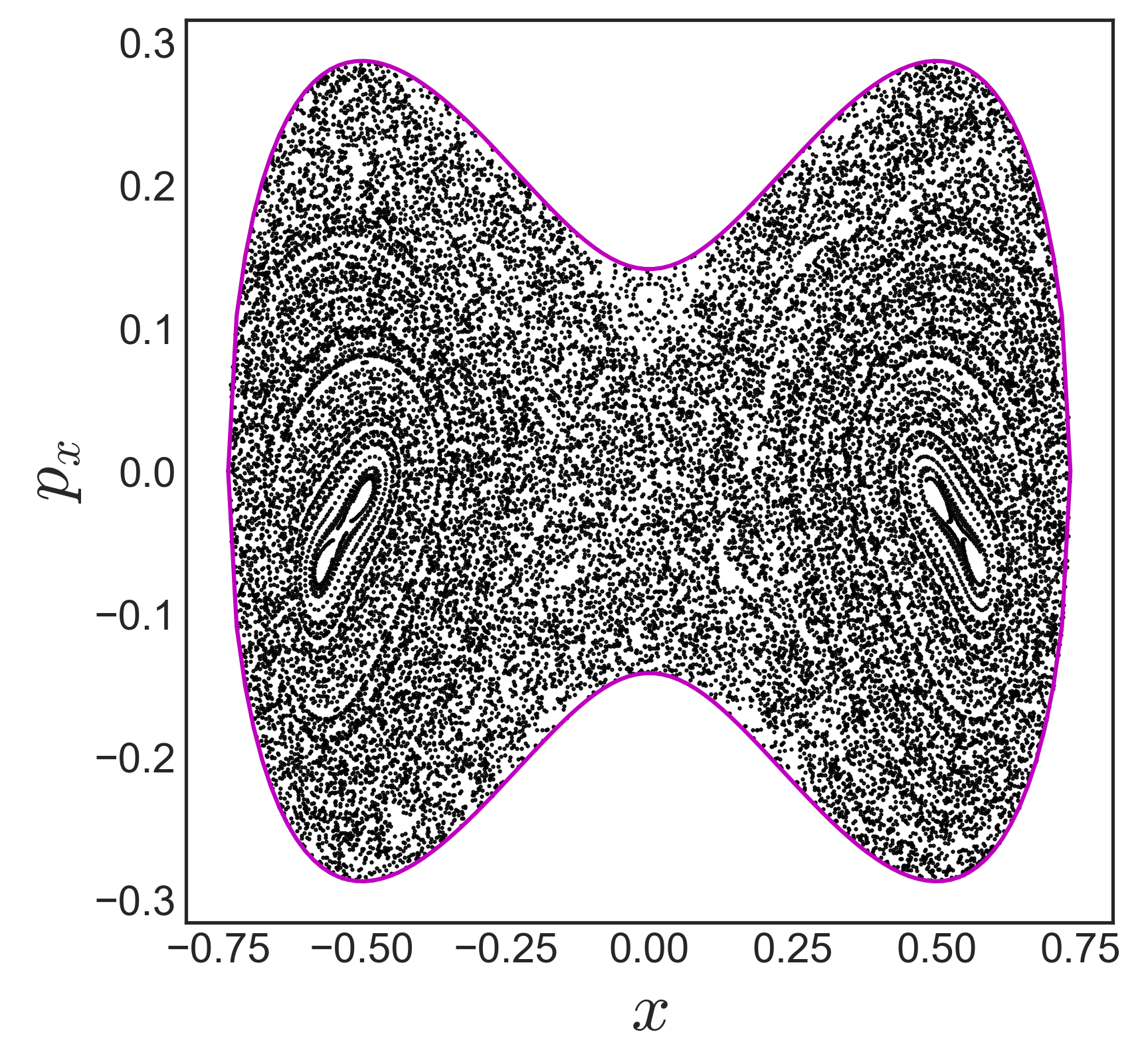}}
	\subfigure[$(1, 1, 4, 0.1)$]{\includegraphics[width=0.48\textwidth]{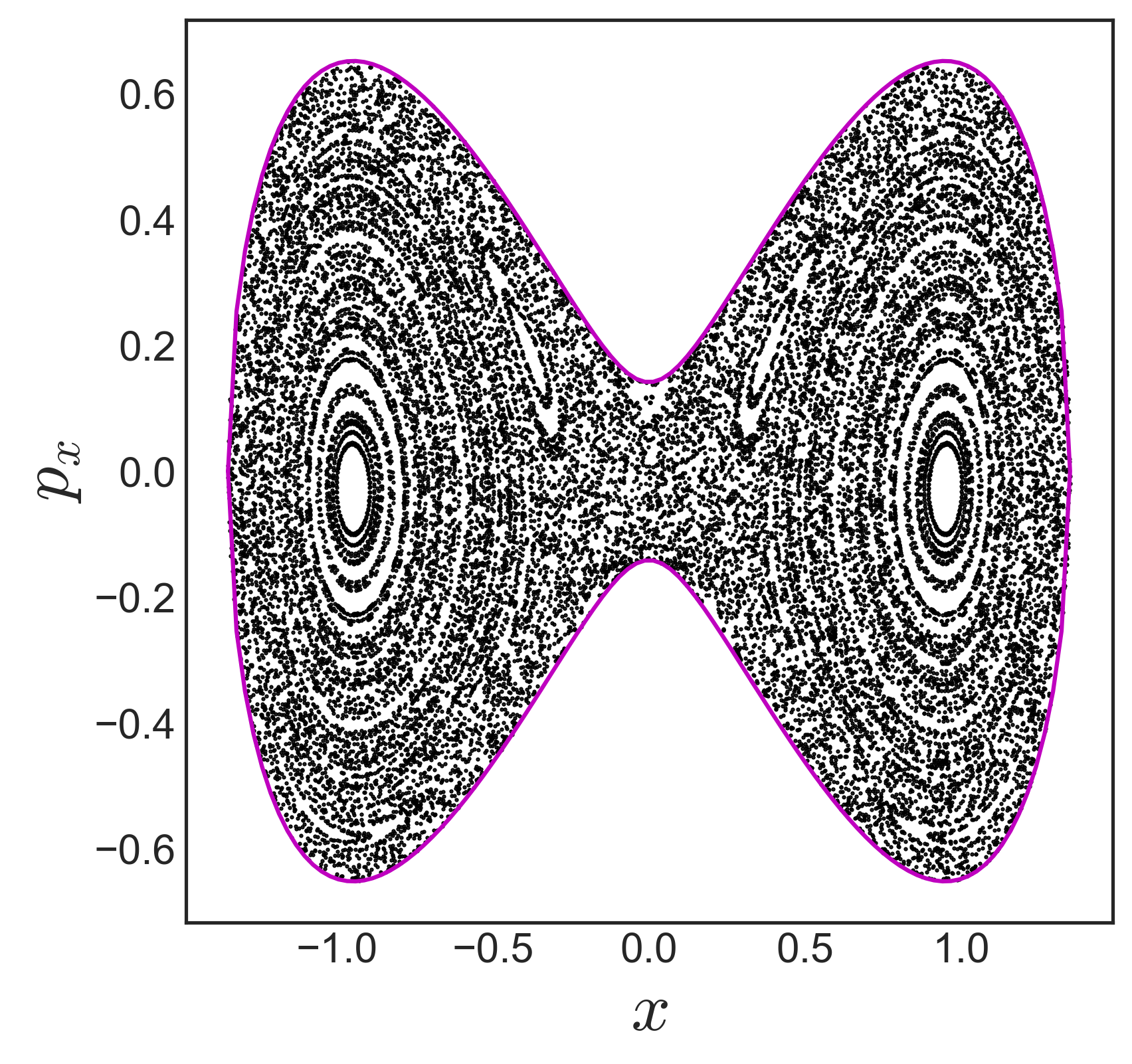}}
	\subfigure[$(1, 1, 0.5, 0.5)$]{\includegraphics[width=0.48\textwidth]{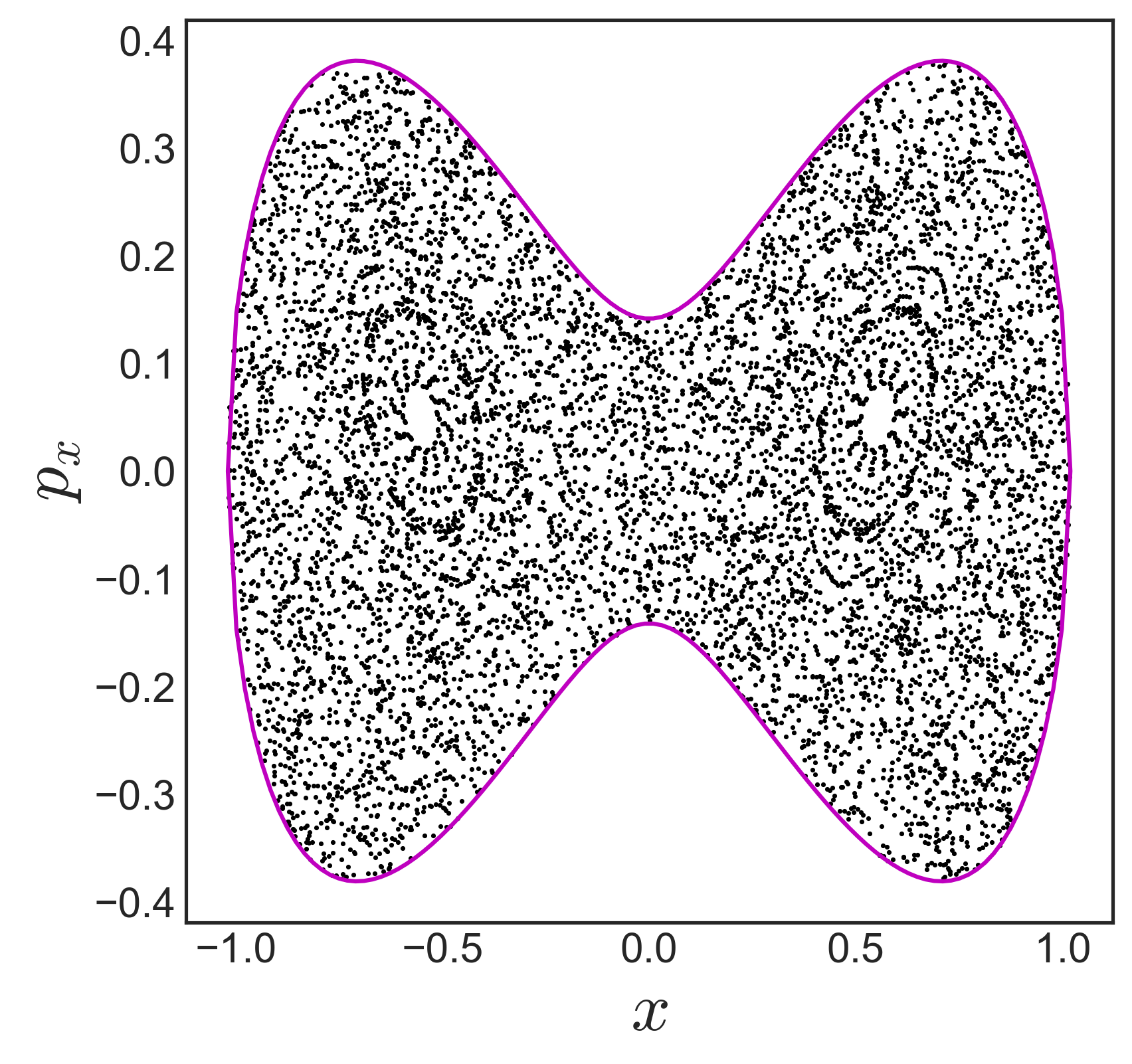}}
	\caption{Poincar\'e surface of sections for the double-well system with different system, bath, and interaction coefficients to illustrate the coexistence of regular and chaotic trajectories in an isomerizing system. All the systems are at same excess energy, $\Delta E = 0.01$, $m_s = m_b = 1.0$ with $(\alpha, \beta, \omega, \epsilon)$ value shown as tuples in the caption for the surface defined by Eqn.~\eqref{eqn:doublewell_xpx_sos}. We note that even though the total energy is constant in all the cases shown here the bottleneck width for these parameters is altered because of the changes in the potential energy surface.}
	\label{fig:psect_bilinear}
\end{figure}

\textbf{Poincar\'e sections.} 
A preliminary step in studying the global dynamics of the double well system shown in Fig.~\ref{fig:pes_contourf_bilinear}(b) is to compute the Poincar\'e section of trajectories by defining the two dimensional surface
\begin{equation}
	\mathcal{U}^+_{xp_x} = \left\{ (x,y,p_x,p_y) \in \mathbb{R}^4 \; | \; y = 0, \, \dot{y}(x,y,p_x; E) > 0 \right\}
	\label{eqn:doublewell_xpx_sos}
\end{equation}
and the resulting Poincar\'e surface of sections are shown in Fig.~\ref{fig:psect_bilinear}. We can observe a mixture of regular and chaotic dynamics for the illustrative values of parameters and for the same total energy, $E = 0.02$, used in the Poincar\'e sections.
Although, the Poincar\'e sections reveal the structure of the global dynamics, it can not be applied to open potential wells since the trajectories that leave the potential well (either in short or long times) become indistinguishable due to white regions in the Poincar\'e surface of section. In addition, when the dynamics occurs on a surface of dimension greater than two, a two-dimensional Poincar\'e section can not reveal the global dynamics because trajectories may not cross the surface of section. In this case, a codimension-1 Poincar\'e section becomes impossible to visualize completely without taking further sections or projections. For these settings of open potential wells and high-dimensional phase space, another trajectory diagnostic is to be relied upon and related results and discussions can be found in Refs.~\cite{junginger2017chemical,craven2017lagrangian,naik2019bfinding,naik2020detecting}.

\textbf{Phase space dividing surface and invariant manifolds associated with the unstable periodic orbit.}
For the uncoupled system, $\varepsilon = 0$, the bifurcation condition becomes $\alpha = 0$. In addition, the uncoupled system is integrable and we can derive an explicit analytical formula for the phase space dividing surface (PSDS). We note here that the dividing surface constructed in the phase space has the locally no-recrossing property~\cite{waalkens2004direct} and in general, trajectories will show global recrossings of the PODS due to the Poincar{\'e} recurrence theorem~\cite{wiggins2003applied} when the energy surface is bounded.
\begin{enumerate}
	\item \textbf{Case II and IV}: When $\alpha > 0$ and $\beta \gtrless 0 $, the dividing surface anchored at the UPO (which is a normally hyperbolic invariant manifold (NHIM)~\cite{wig2016} in a two DOF system) is defined by the condition $x = 0, p_x = 0$ on the energy surface. The dividing surface which is a 2-sphere ($\mathbb{S}^2$) has forward ($p_x(y,p_y;e) < 0$) and backward ($p_x(y,p_y;e) > 0$) hemispheres given by
	\begin{align}
		\text{forward/backward DS:} \quad & \left\{ (x,y,p_x,p_y) \in \mathbb{R}^4 \, | \, x = 0, p_x(y,p_y;e) \lessgtr 0 \right\} \label{eqn:forwbackds_quartic_uncoupled_alp>0}
	\end{align}
	where $e$ is the total energy of the system and $p_x(y,p_y;e)  = \pm \sqrt{ 2e - p_y^2  -\omega y^2}$ is calculated from the fixed energy constraint. The forward reaction occurs when the trajectory crosses the forward DS and the backward reaction occurs when the trajectory crosses the backward DS. In case IV, which is simply a two DOF saddle, the energy surface is unbounded (open potential well), trajectories go-off to infinity and do not return to recross the PODS. However, trajectories will recross the PODS in case II due to Poincar\'e recurrence theorem and we will characterize this using a quantitative measure of this recrossing called \emph{gap times}~\cite{ezra_microcanonical_2009}. We note here that the global recrossing property is independent of the existence of chaos in a dynamical system. The uncoupled system is an example of that since it is integrable, and hence can not exhibit chaotic dynamics. 
	
	
	\item \textbf{Case III}: When $\alpha < 0$ and $\beta <0$, there are two index-1 saddles, and hence two corresponding UPOs ($\Gamma_{1,2}$) and two PSDSs. The UPOs are defined by the condition $x = \pm \sqrt{\alpha/\beta}, p_x = 0$ on the energy surface and the PSDSs are given by: 
	\begin{align}
		\text{ forward/backward $\rm{DS}_-$:} \quad & \left\{ (x,y,p_x,p_y) \in \mathbb{R}^4 \, | \, x = \sqrt{\alpha/\beta}, p_x(y,p_y;e) \lessgtr 0 \right\} \label{eqn:forwbackds1_uncoupled_alp<0}\\
		\text{ forward/backward $\rm{DS}_+$:} \quad & \left\{ (x,y,p_x,p_y) \in \mathbb{R}^4 \, | \, x = -\sqrt{\alpha/\beta}, p_x(y,p_y;e) \lessgtr 0 \right\} \label{eqn:forwbackds2_uncoupled_alp<0}
	\end{align}
	where $e$ is the total energy of the system and $p_x(y,p_y;e)  = \pm \sqrt{ 2e - p_y^2  -\omega y^2}$ is calculated from the fixed energy constraint. 
\end{enumerate}

We note that the \textbf{case I} $\alpha < 0, \beta > 0$ only has one equilibrium point of centre-centre stability, so the system does not show transition (or reaction) dynamics as defined in for this model Hamiltonian.

\begin{figure}
	\centering
	\subfigure[]{\includegraphics[width=0.45\textwidth]{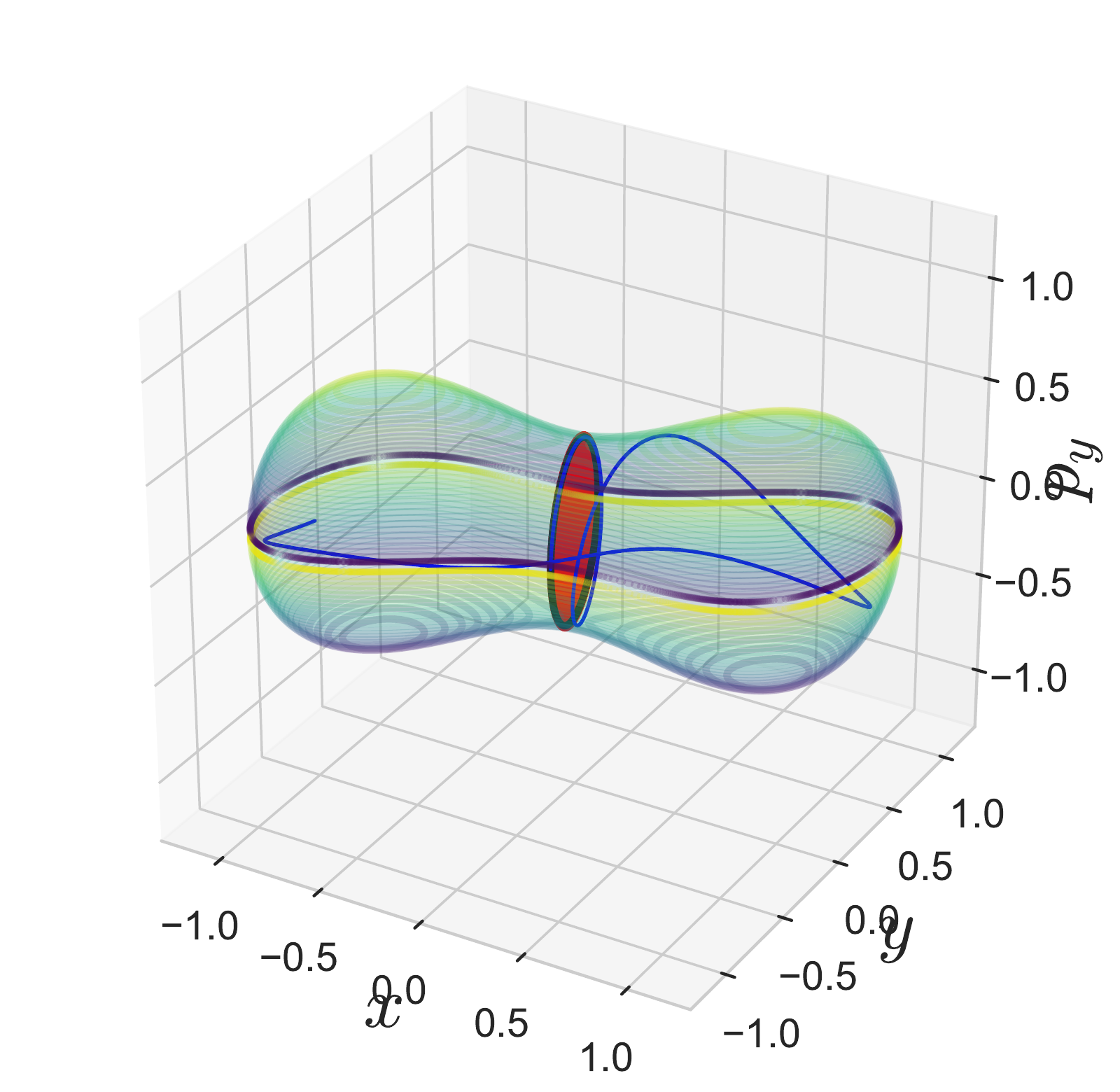}}
	\subfigure[]{\includegraphics[width=0.45\textwidth]{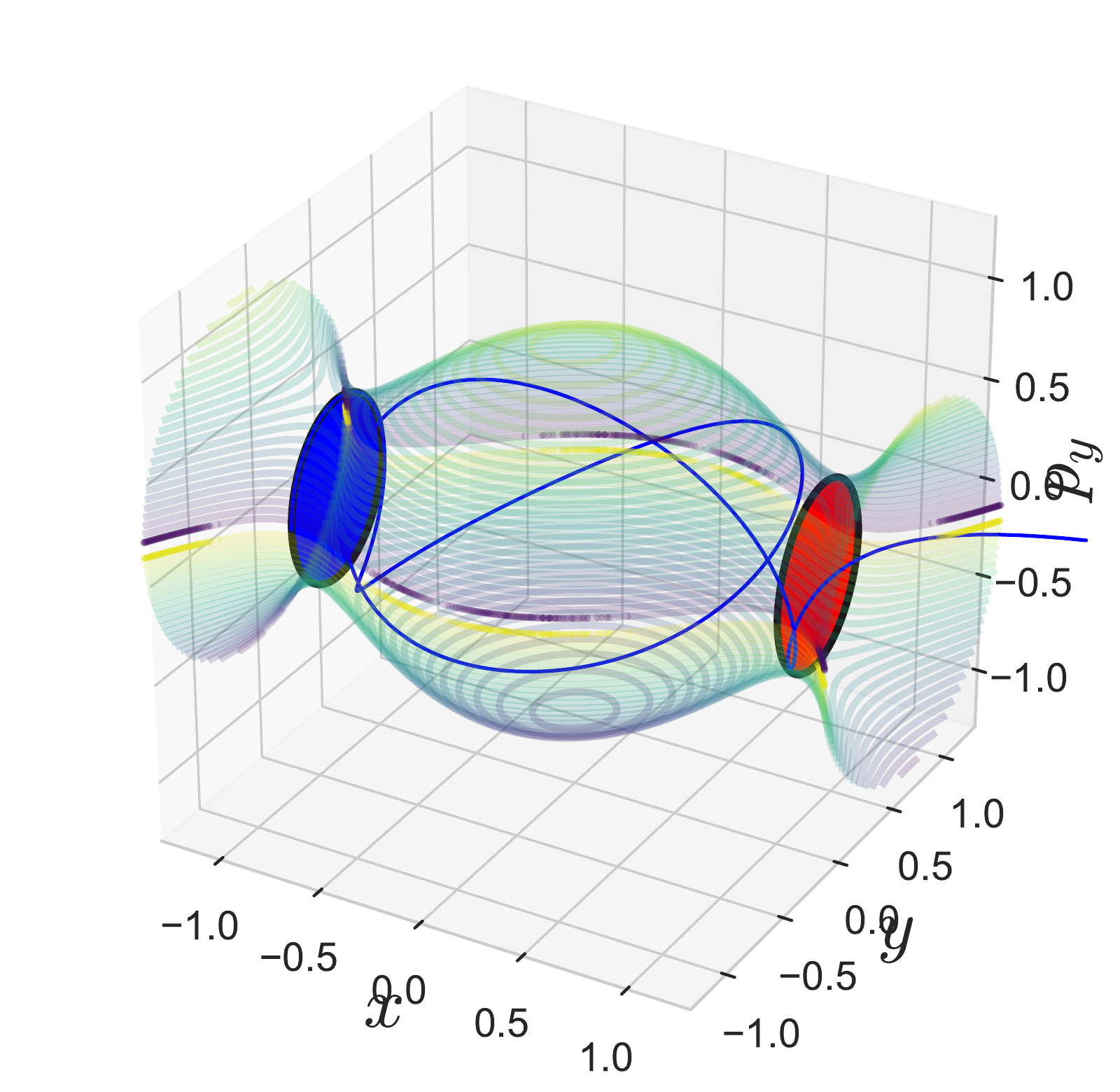}}
	\caption{\textbf{DS}(blue, red disks), NHIM(black) on the energy surface for the coupled system with (a)  $\alpha=\beta=1,\varepsilon=0.5$ (b) negative $\alpha=\beta=-1,\varepsilon=0.2$. Other parameters are chosen as $\Delta E =0.1,\omega=1$.}
	\label{fig:DS_visual_bilinear}
\end{figure}

For the coupled system, the phase space dividing surface is sampled using the algorithm in Ref.~\cite{Ezra_sampling2018} after obtaining the unstable periodic orbit for a given energy, $\Gamma(E)$ using the open-source python package, uposham~\cite{Lyu2020}. Let us denote the unstable periodic orbit by $\Gamma(E) := (\overline{x},\overline{y}, \overline{p}_x, \overline{p}_y) \in \mathbb{R}^4$ or by $\Gamma_{-,+}$ when there are two unstable periodic orbits at $x < 0$ and $x > 0$. 
\begin{enumerate}
	\item \textbf{Case II and IV}: When $\alpha > \omega \varepsilon/(\omega + \varepsilon)$, the PODS is given by
	\begin{align}
		\text{ forward/backward DS:} & \; \left\{ (x,y,p_x,p_y) \in \mathbb{R}^4 \, \vert \, p_x(\overline{x},\overline{y},\overline{p_y};e) \lessgtr 0 \right\} \label{eqn:forwbackds_quartic_coupled_alp>crtical}
	\end{align}
	where $p_x  = \pm \sqrt{ 2(e - p_y^2-V(\overline{x},\overline{y}))}$ is the fixed energy constraint and $\overline{x},\overline{y} \in \mathrm{UPO}$. 
	
	\item \textbf{Case III}: When $\alpha < \omega \varepsilon/(\omega + \varepsilon), \beta <0$, the PSDSs are given by
	\begin{align}
		\text{ forward/backward $\rm{DS}_-$:} & \; \left\{ (x,y,p_x,p_y) \in \mathbb{R}^4 \, \vert \, p_x(\overline{x}_-,\overline{y}_-,p_y;e) \lessgtr 0 \right\} \label{eqn:forwbackds1_quartic_coupled_alp<crtical}\\
		\text{ forward/backward $\rm{DS}_+$:} & \; \left\{ (x,y,p_x,p_y) \in \mathbb{R}^4 \, \vert \, p_x(\overline{x}_+,\overline{y}_+,p_y;e) \lessgtr 0 \right\} \label{eqn:forwbackds2_quartic_coupled_alp<crtical}
	\end{align}
	where $p_x  = \pm \sqrt{ 2(e - p_y^2-V(\overline{x}_{-,+},\overline{y}_{-,+}))}$ is the fixed energy constraint and $\overline{x}_{-,+},\overline{y}_{-,+} \in \Gamma_{-,+}$. 
\end{enumerate}
The phase space dividing surfaces on the energy surface for the coupled systems are shown in Fig.~\ref{fig:DS_visual_bilinear} . We note that the PODS have a geometry of discs in the two configuration space coordinates with one momentum coordinate. 

\begin{figure}[!ht]
	\centering
	\subfigure[]{\includegraphics[width=0.52\textwidth]{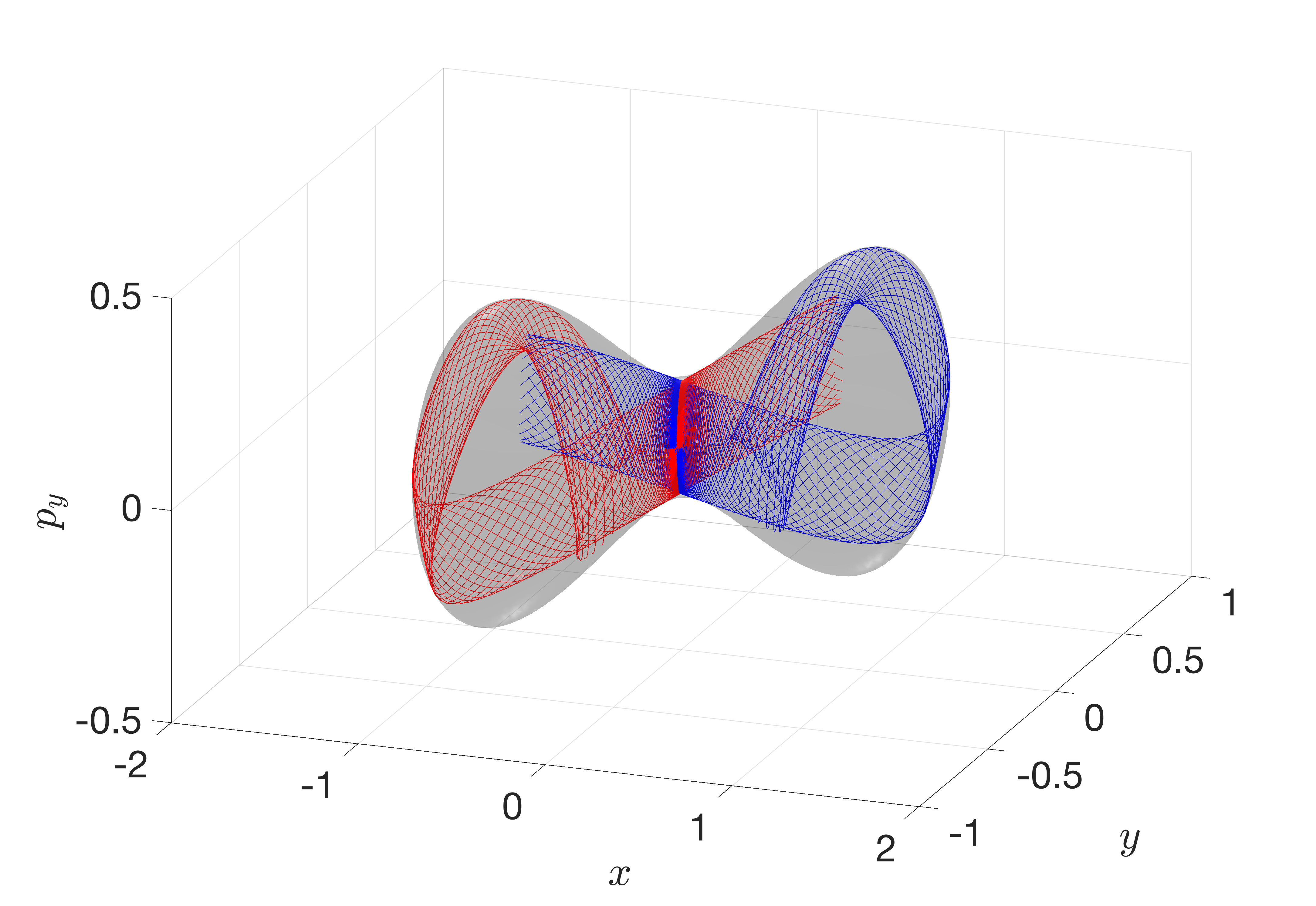}}\\
	\subfigure[]{\includegraphics[width=0.50\textwidth]{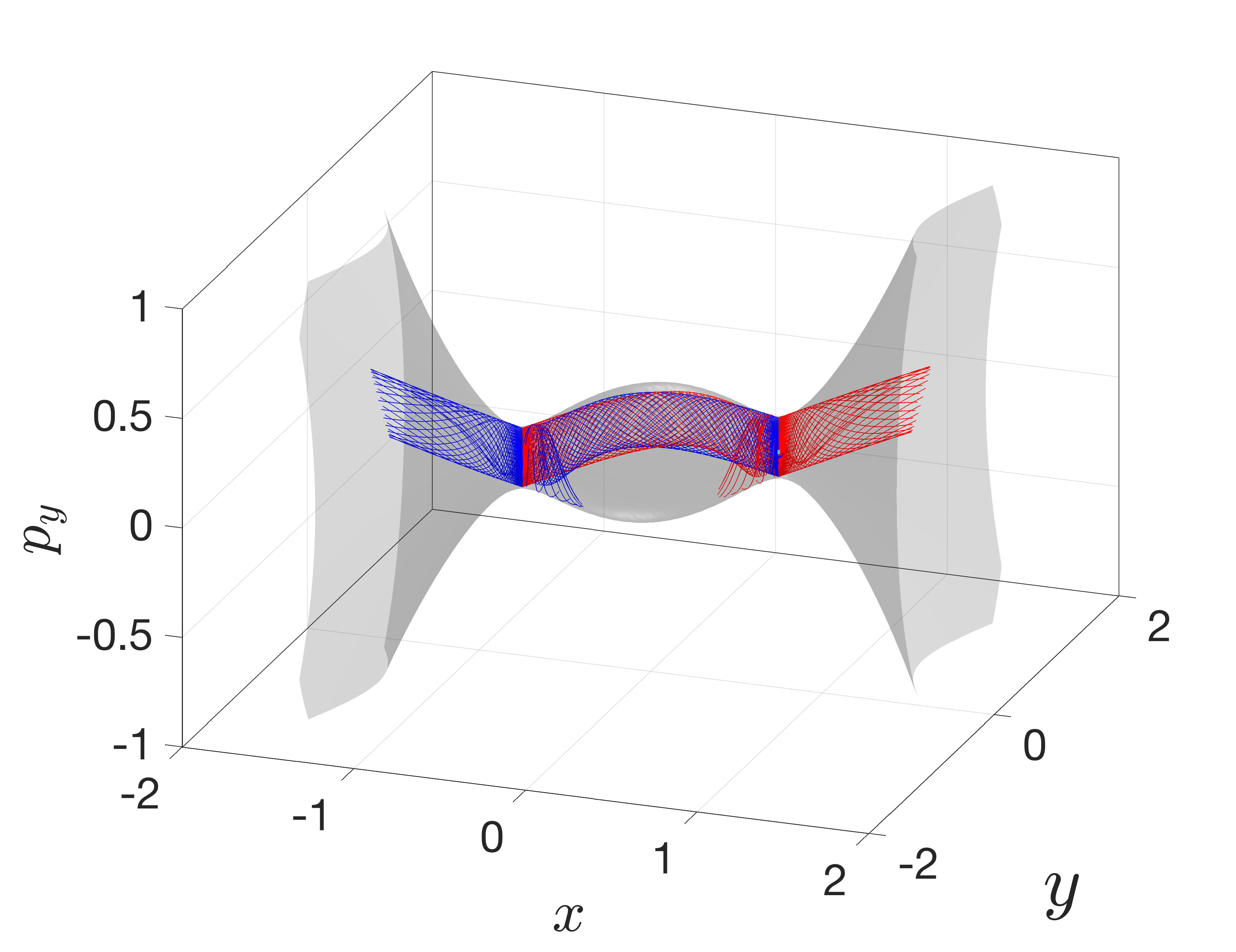}}
	\subfigure[]{\includegraphics[width=0.45\textwidth]{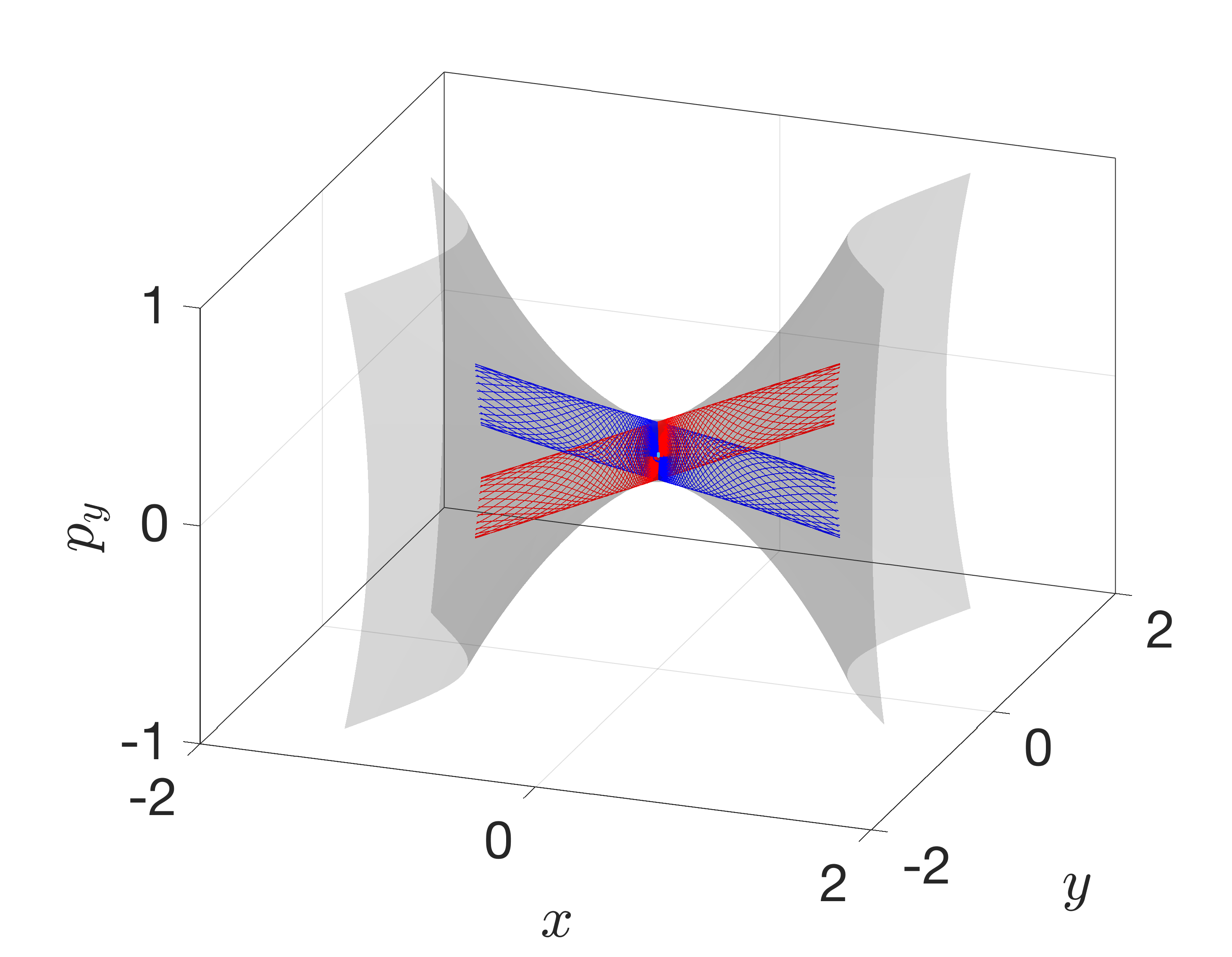}}
	\caption{\textbf{Cylindrical invariant manifolds and energy surface} for the three reactive systems due to the pitchfork bifurcation. Stable (blue) and unstable (red) manifolds of the unstable periodic orbit associated with the index-1 saddle (NHIM for two DOF systems) in the bottleneck while the energy surface is shown as the grey surface. (a) \textbf{Case II:} $\alpha = 1.0, \beta = 1.0$ (b) \textbf{Case III:} $\alpha = 0.1, \beta = -1.0$ (c) \textbf{Case IV:} $\alpha = 1.0, \beta = -1.0$. Other parameters, $\omega = \varepsilon = 1.0$, and excess energy, $\Delta E = 0.01$, are same for all cases.}
	\label{fig:invariant_manifolds_energysurface_E1e-2_bilinear}
\end{figure}

From a phase space perspective of transition dynamics across index-1 saddles, the phase space dividing surface can be viewed as a local barrier between reactive and non-reactive trajectories. However, the global barriers between the two type of trajectories are the invariant (stable and unstable) manifolds associated with the UPO. These stable and unstable invariant manifolds have a cylindrical geometry, that is $\mathbb{R}^1 \times \mathbb{S}$, and form impenetrable barriers between the reactive and non-reactive trajectories and shown in Fig.~\ref{fig:invariant_manifolds_energysurface_E1e-2_bilinear}.


\subsection{Gap time distributions of the transition}

We use the definition of the gap time which is defined as the time interval between two successive crossings of the phase space DS~\cite{ezra_microcanonical_2009}. For \textbf{case II}, that is when $\beta > 0, \alpha > (\omega \varepsilon)/(\omega+ \varepsilon)$, the gap time distribution is calculated by launching trajectories from the DS with $p_x > 0$ such that the trajectory enters the reactants region. Then, the time interval until the next crossing of the DS represents the gap time or the first passage time. We note that local recrossings are impossible since the initial conditions are on a phase space dividing surface. For this case with double-well and a bounded energy surface, when a trajectory recrosses the DS and enters the products region, the $x-$momentum changes sign and $x < 0$. This happens when the trajectory crosses the forward DS given by $p_x < 0$ in the Eqn.~\eqref{eqn:forwbackds_quartic_coupled_alp>crtical}. We perform this calculation for an ensemble of initial conditions at a fixed energy and we record the gap time for all the trajectories starting on the PSDS with positive $p_x$ within 100 time units (same for all the case discussed here), and call this the microcanonical gap time distribution~\cite{ezra_microcanonical_2009}.

For \textbf{case III}, that is when $\beta < 0, \alpha < (\omega \varepsilon)/(\omega+ \varepsilon)$, there are two phase space dividing surfaces associated with the two index-1 saddles located at $x < 0$ and $x > 0$. Thus, for each PSDS there are two gap time distributions: one for crossing the PSDS where the ensemble gets initialized and one for crossing the other PSDS. A comparison of these two gap times is beyond the scope of this article but interested readers can see an examples of a similar analysis for a three degrees of freedom model Hamiltonian of HCN isomerization in Ref.~\cite{waalkens_phase_2004}. We note here that out of the two phase space dividing surfaces, we choose $\rm{DS}_-$ to initialize an ensemble at a fixed energy with $p_x > 0$ and integrate until they either cross $\rm{DS}_-$ or $\rm{DS}_+$. We record only the time to cross the PSDS $\rm{DS}_+$ within 100 time units.

\begin{figure}[!ht]
	\centering
	\subfigure[]{\includegraphics[width=0.45\textwidth]{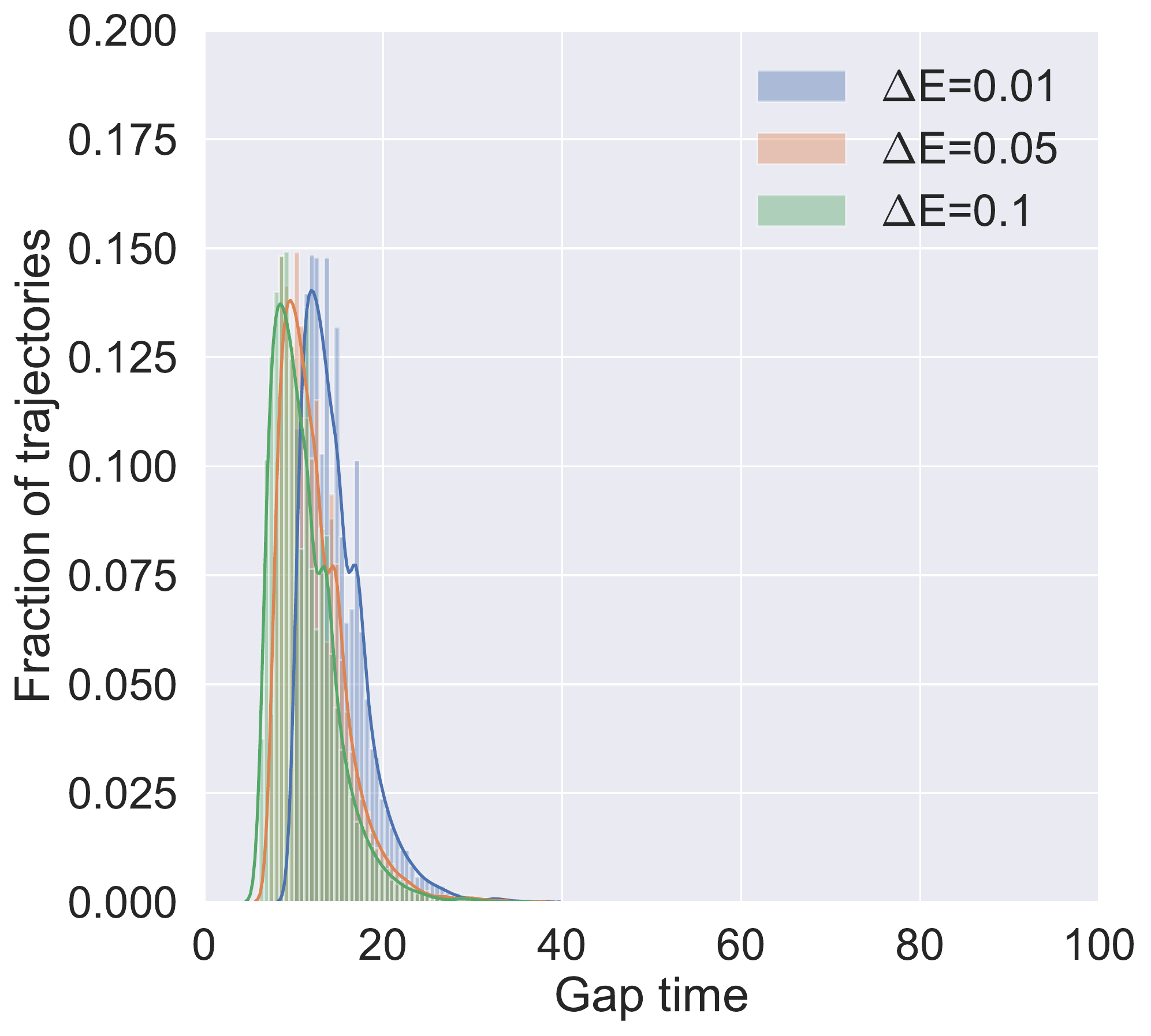}}
	\subfigure[]{\includegraphics[width=0.45\textwidth]{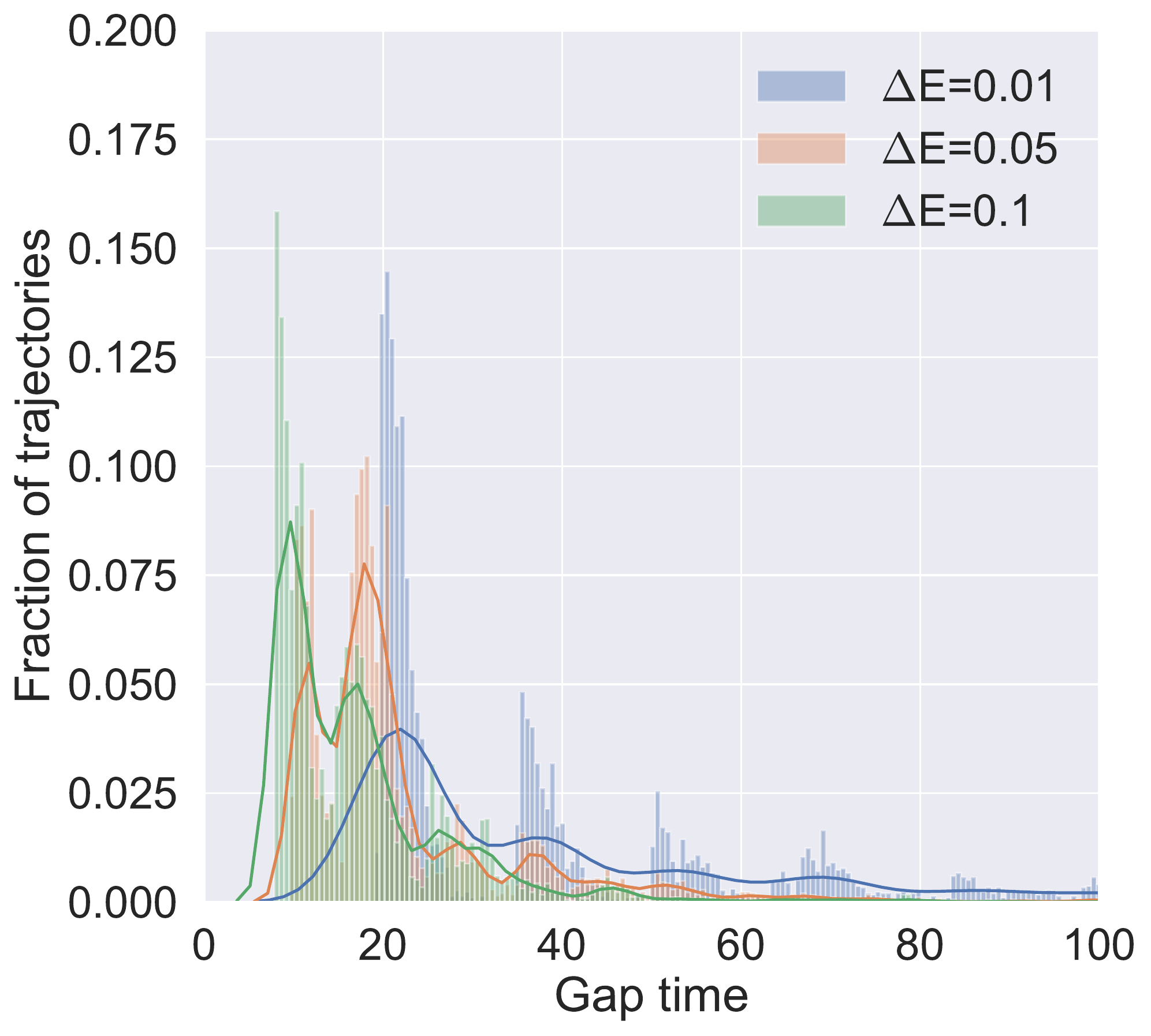}}
	\caption{\textbf{Gap time distributions for case II} for the (a) uncoupled $\epsilon=0$ and the (b) coupled $\epsilon=0.5$ system. Parameters in the potential energy are chosen as $\alpha=1, \beta=1, \omega=1$.}
	\label{fig:gap_time_alp1_beta1_omega1_bilinear}
\end{figure}
\begin{figure}[!ht]
	\centering
	\subfigure[]{\includegraphics[width=0.45\textwidth]{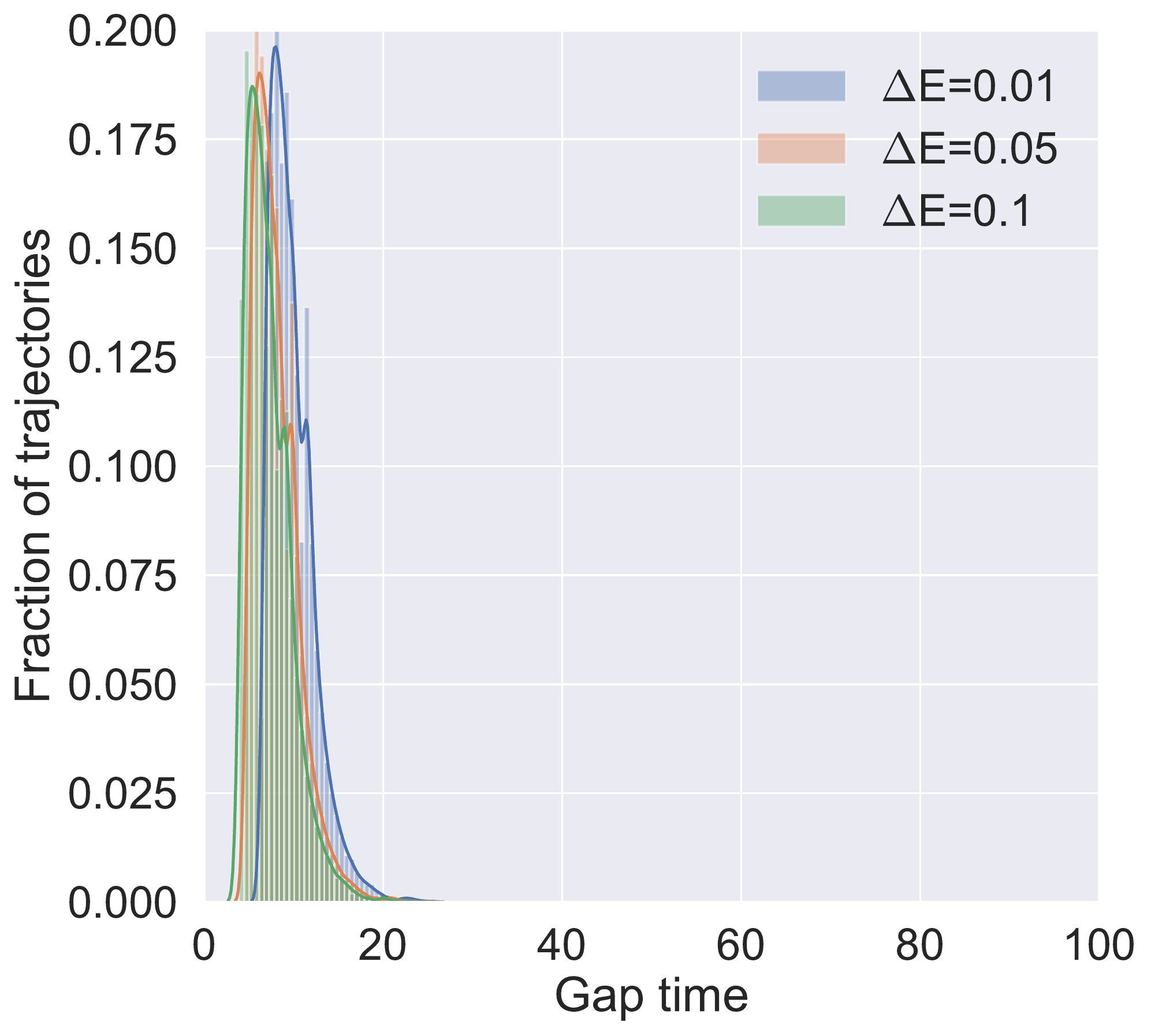}}
	\subfigure[]{\includegraphics[width=0.45\textwidth]{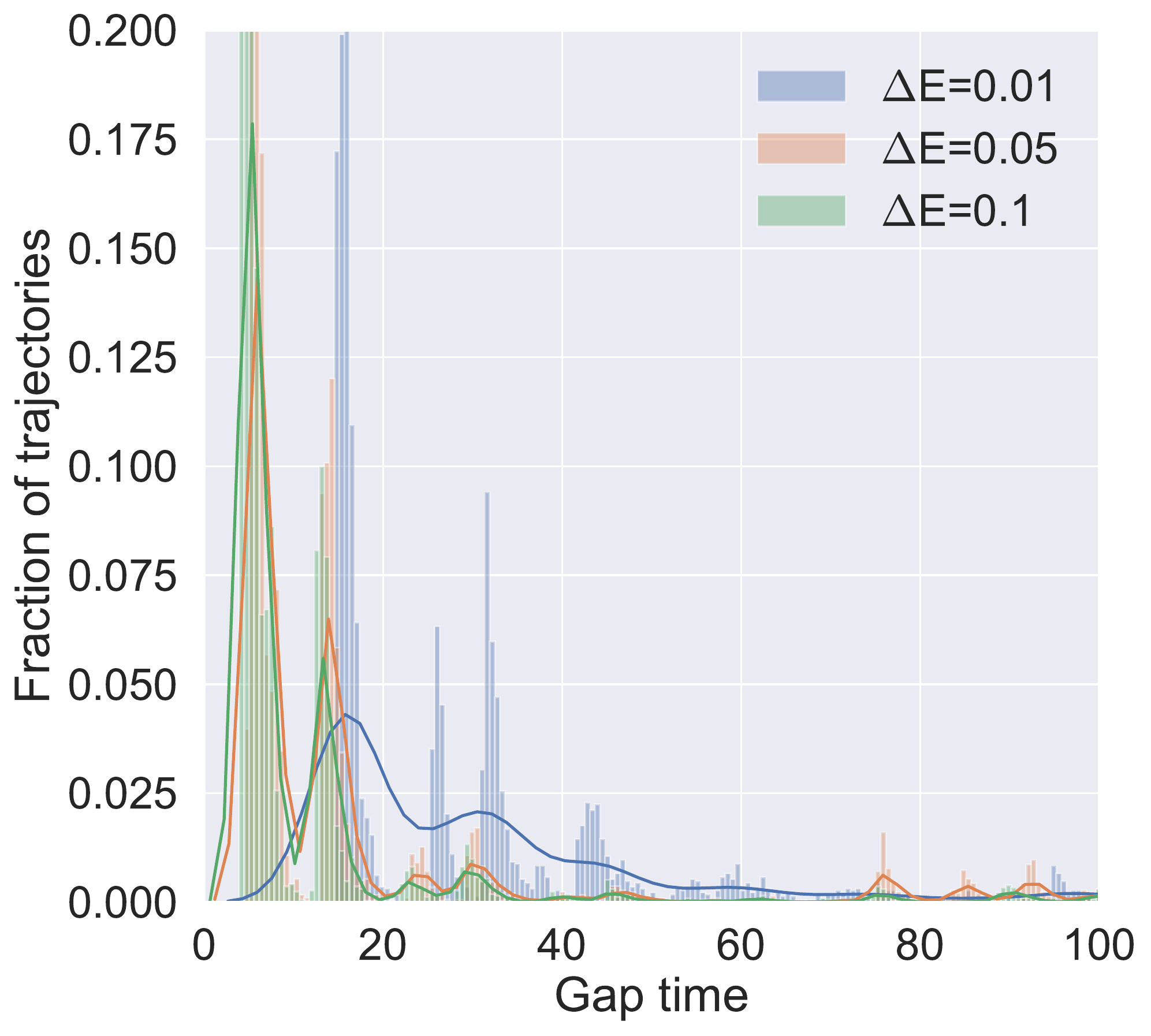}}
	\caption{\textbf{Gap time distributions for case III} for the  (a) uncoupled $\epsilon=0$ (b) and the coupled $\epsilon=0.4$ system. Parameters in the potential energy are chosen as $\alpha=-1, \beta=-1, \omega=1$.}
	\label{fig:gap_time_alpminus1_betaminus1_omega1_bilinear}
\end{figure}

We use the algorithm in Ref.~\cite{Ezra_sampling2018} to sample points on the phase space DS (which is a hemisphere of the 2-sphere, $\mathbb{S}^2$) to calculate the gap time distribution. \emph{Firstly}, we obtain the UPO associated with the index-1 saddle (or saddles in case III) equilibrium point for a given energy $e$.  \emph{Second}, we use the configuration space projection of the UPO $\{(\overline{x}_i,\overline{y}_i)\}$. For each point on the configuration space $\{(\overline{x}_i,\overline{y}_i)\}$, the $(p_x, p_y)$ satisfies the circle given by
\begin{equation}
	\frac{1}{2}(p_x^2+p_y^2) = e - V(x_i,y_i).
\end{equation}
Therefore, the maximum value of $p_x$ is
\begin{equation}
	p_x^{max} = \sqrt{2(e - V(x_i,y_i))}
\end{equation}
and we select $p_{x,i}$ uniformly from the interval $[-p_x^{max},p_x^{max}]$. Using the definition of the Hamiltonian, we calculate the value of $p_{y,i}$, which is either positive or negative. We note that $p_y^{max} = p_x^{max},\ p_y \in \left[-p_y^{max},p_y^{max}\right]$. \emph{Thirdly}, we generate a set of microcanonical ensemble $\{( x_i,y_i,p_{x,i},p_{y,i})\}_e$ on the phase space DS at a total energy $e$. We select the initial conditions with positive $p_x$ and integrate so that the trajectories enter the reactants ($x > 0 $ in case II or between the index-1 saddles in case III) region, and leave the reactants region by either crossing the same phase space DS (as in case II) or the other phase space DS (as in case III). \emph{Finally,} the time interval between entering and exiting the reactants region is the gap time of the initial condition. For the coupled systems, the UPO is computed using an open-source python library for computing unstable periodic orbits in two DOF Hamiltonian systems~\cite{Lyu2020}. Thus, for case II, to simplify this detection of crossing the forward DS, we use an alternate boundary at $x = -1$ line which ensures that the trajectory has left the reactants region and is in the products region. The gap time distributions for $\Delta E = 0.01, 0.05, 0.1$ are shown in Fig.~\ref{fig:gap_time_alp1_beta1_omega1_bilinear}(a) and~\ref{fig:gap_time_alp1_beta1_omega1_bilinear}(b) for the uncoupled and coupled systems, respectively. For the case III, we use an alternate boundary at $x = 1.5$ which ensures that the trajectory has left the reactants region and is in the products region. The gap time distributions for $\Delta E = 0.01, 0.05, 0.1$ are shown in Fig.~\ref{fig:gap_time_alpminus1_betaminus1_omega1_bilinear}(a) and~\ref{fig:gap_time_alpminus1_betaminus1_omega1_bilinear}(b) for the uncoupled and coupled systems, respectively.

We observe some common features in the gap time distributions for both the cases. We see from the figures that the gap time distributions have a unimodal shape for the uncoupled system except for a time shift of the peak of the distribution. For the coupled systems, the gap time distributions have multiple pulses and the time of the first pulse appears at time instant that decreases with increasing excess energy. Thus, as $\Delta E$ increases, most trajectories cross the phase space DS to the products region at short timescales compared to the timescale in which all the trajectories exit the reactant region. This is observed across systems with single and multiple phase space dividing surfaces~\cite{ezra_microcanonical_2009,waalkens_phase_2004}.

\subsection{Directional flux of the transition}
The directional flux of transition trajectories through the phase space dividing surface is given by the action of the normally hyperbolic invariant manifold (NHIM) which is of dimension $\mathbb{S}^{2N - 3}$ in a $N$ degrees of freedom~\cite{waalkens2004direct,Waalkens2005}. In the two DOF system, the NHIM is an unstable periodic orbit (UPO) and the directional flux given by the action simplifies to the line integral:
\begin{equation}
	Q = \int_{\rm{UPO}} \mathbf{p} \cdot d\mathbf{q}.
\end{equation}
For the uncoupled systems, the UPO can be expressed implicitly by substituting the condition $x = 0, p_x = 0$ in the Hamiltonian.
\begin{enumerate}
	\item \textbf{Case II and IV:} 
	\begin{align}
		\text{ UPO:} \; \left\{ (x,y,p_x,p_y) \in \mathbb{R}^4 \, | \, x = 0, p_x = 0, \dfrac{1}{2} p_y^2 +\dfrac{\omega}{2} y^2 = e \right\} 
	\end{align}
	Thus, the directional flux, $Q$, is given by
	\begin{align}   
		Q & = \int_{\text{UPO}} \mathbf{p} \cdot d\mathbf{q} = \int_{0}^{T} \mathbf{p} \cdot \dfrac{d\mathbf{q}}{dt} \ dt = \int_{0}^{T} (p_x \dot{x} + p_y \dot{y} ) \ dt = \int_{0}^{T}  p_y \dot{y} \ dt  = \int_{0}^{T}  p_y^2 \ dt = Te = \frac{2 \pi e}{\sqrt{\omega}} \label{eqn:caseII_IV_flux}
	\end{align}    
	where $T$ is the time period of the UPO.
	\item \textbf{Case III:} When $\alpha<0$ and $\beta<0$, the directional flux through either DS, is given by
	\begin{align}   
			Q & = \int_{\text{UPO}} \mathbf{p} \cdot d\mathbf{q} = \int_{0}^{T} \mathbf{p} \cdot \dfrac{d\mathbf{q}}{dt} \ dt = \int_{0}^{T} (p_x \dot{x} + p_y \dot{y} ) \ dt = \int_{0}^{T}  p_y \dot{y} \ dt = \int_{0}^{T}  p_y^2 \ dt = T \Delta E  = \frac{2 \pi \Delta E}{\sqrt{\omega}} \label{eqn:caseIII_flux}
	\end{align}    
	where $T$ is the time period of the UPO. 
\end{enumerate}
 In Eqns.~\eqref{eqn:caseII_IV_flux},~\eqref{eqn:caseIII_flux}, $p_x \dot{x}$ vanishes because the vector field is tangential on the UPO. We calculate the integral by expressing $p_y$ in terms of $t$, which can be solved explicitly for the uncoupled Hamilton's equations of motion~\eqref{hameq_2dof}. We note that this expression is the $N = 2$ case for the flux formula in Ref.~\cite{waalkens2004direct}.

We note the directional flux is $\frac{2 \pi \Delta E}{\sqrt{\omega}}$ all three uncoupled systems. For Case II and IV, $\Delta E = e - \mathcal{H}(\mathbf{x}_1^e) = e$ where $e$ is the total energy of the system and $\mathcal{H}(\mathbf{x}_1^e)=0$ is the energy of the saddle equilibrium point. For Case III, $\Delta E = e - \mathcal{H}(\pm \mathbf{x}_2^e) = e + \dfrac{1}{4\beta} \left(\alpha- \frac{\omega\varepsilon}{\omega + \varepsilon} \right)^2$ where $e$ is the total energy of the system and $\mathcal{H}(\pm \mathbf{x}_2^e)=-\dfrac{1}{4\beta} \left(\alpha- \frac{\omega\varepsilon}{\omega + \varepsilon} \right)^2$ is the energy of the saddle equilibrium points.

\textbf{Case II-IV:} For the coupled systems, the directional flux through the phase space DS is computed by evaluating the integral 
\begin{align}   
	Q & = \int_{\text{UPO}} \mathbf{p} \cdot d\mathbf{q} = \int_{0}^{T} \mathbf{p} \dfrac{d\mathbf{q}}{dt} \ dt \\
	Q & = \int_{0}^{T} (p_x \dot{x} + p_y \dot{y} ) \ dt = \int_{0}^{T} (p_x \dot{x} + p_y \dot{y} ) \ dt \\
	Q & =\int_{0}^{T} p_x^2 + p_y^2 \ dt 
\end{align}  


\begin{figure}[!ht]
	\centering
	\subfigure[]{\includegraphics[width=0.45\textwidth]{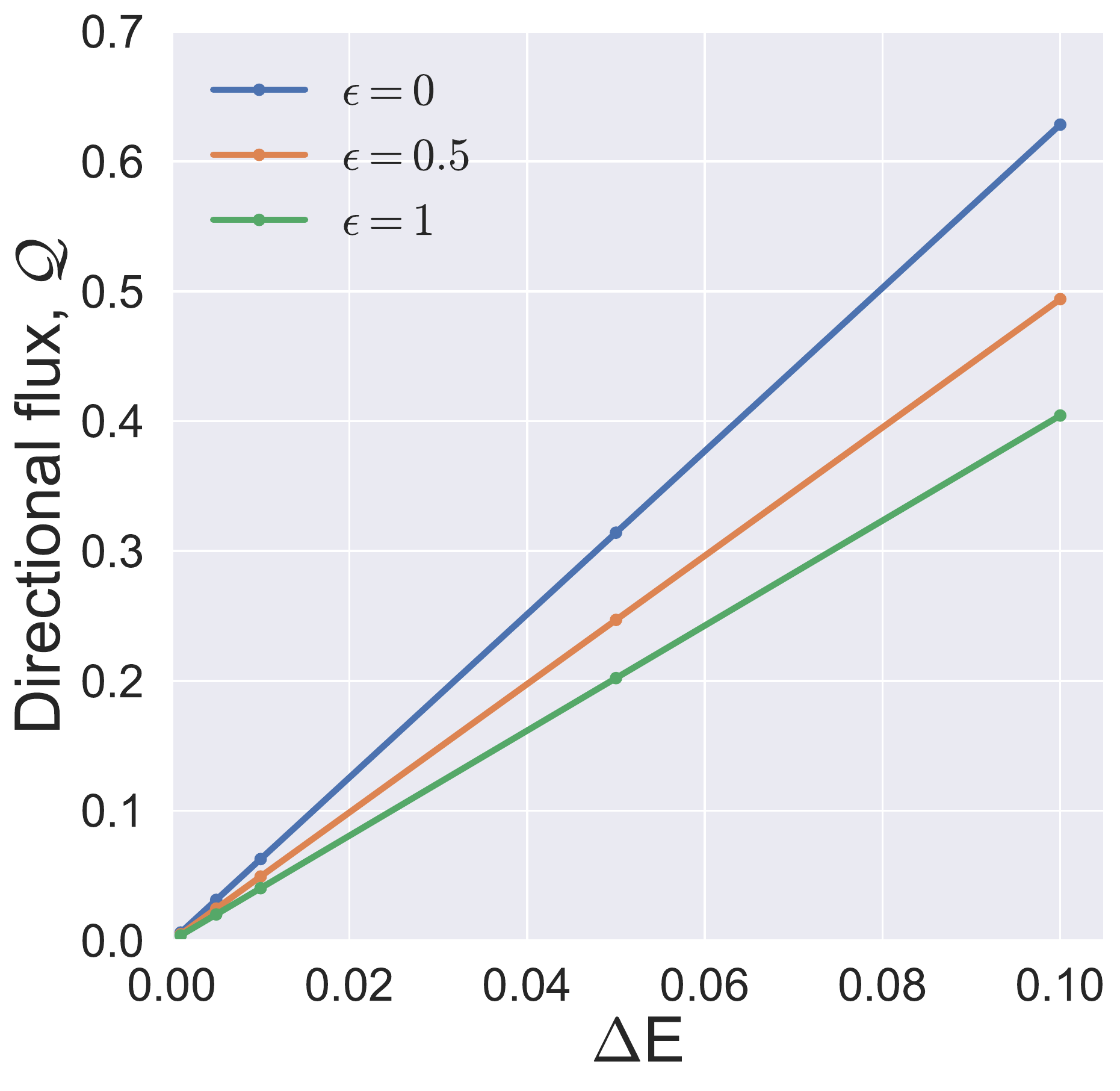}}
	\subfigure[]{\includegraphics[width=0.45\textwidth]{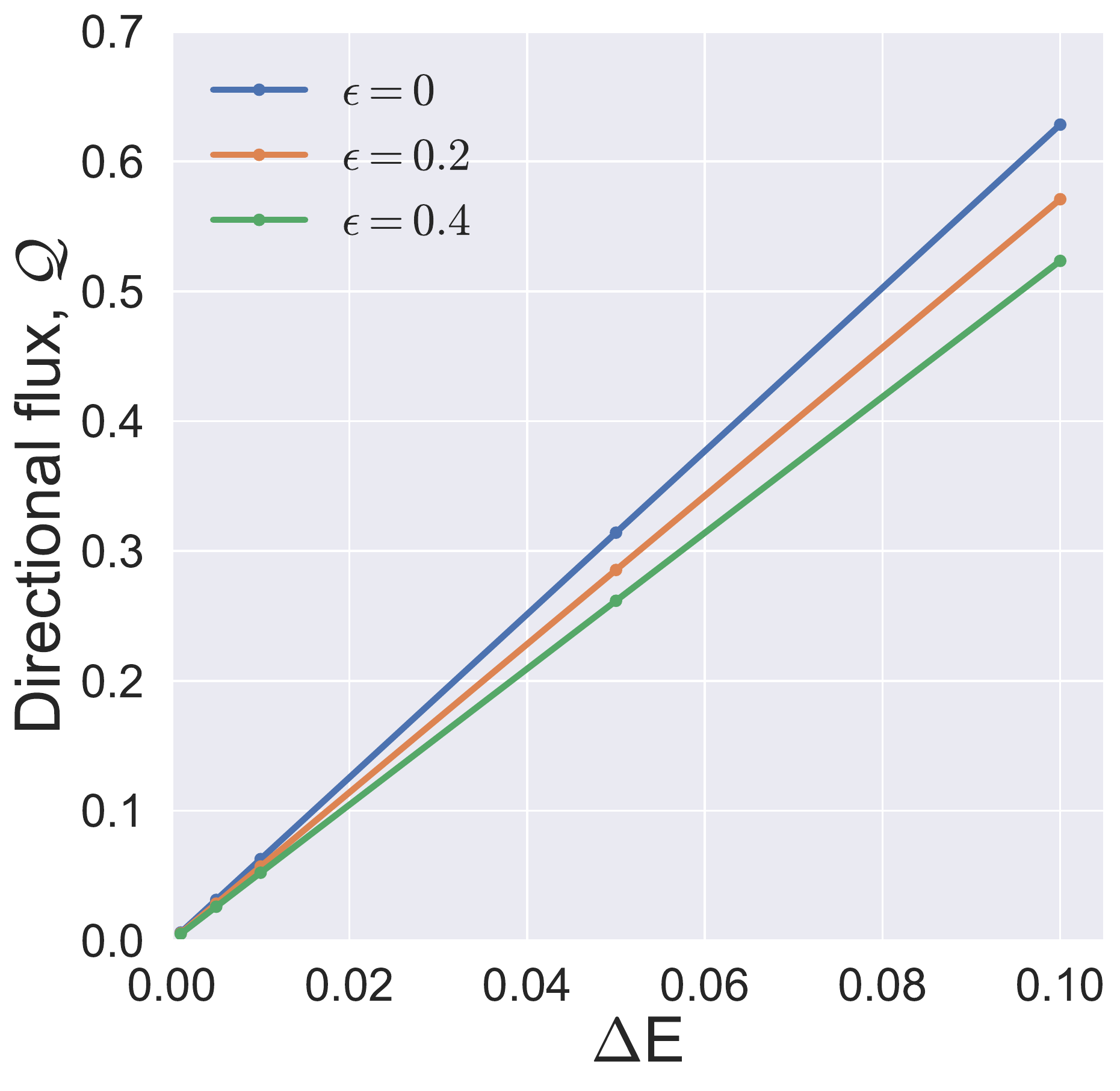}}
	\caption{Directional flux for the transition through the phase space dividing surface for (a) \textbf{case II:} $\alpha=1, \beta=1, \omega=1$, (b) \textbf{case III} and $\rm{DS}_-$: $\alpha=-1, \beta=-1, \omega=1$.} 
	\label{fig:Q_alp1_beta1&alpminus1_betaminus1_deltaE_bilinear}
\end{figure}

In Fig.~\ref{fig:Q_alp1_beta1&alpminus1_betaminus1_deltaE_bilinear}(a) we present the directional flux, $Q$ through the DS for $\alpha=1, \beta=1$ and $\omega=1$ with variation of $\Delta E$. We can see that as $\Delta E$ increases, Q increases linearly. At a fixed value of $\Delta E$, $Q$ decreases as we increase the strength of the coupling $\varepsilon$. In Fig.~\ref{fig:Q_alp1_beta1&alpminus1_betaminus1_deltaE_bilinear}(b) we present the directional flux, $Q$ through the DS for $\alpha=-1, \beta=-1$ and $\omega=1$ with variation of $\Delta E$. We can also see that as $\Delta E$ increases, Q increases linearly. At a fixed value of $\Delta E$, $Q$ decreases as we increase the strength of the coupling $\varepsilon$.

\section{Conclusions and outlook\label{sect:conclusions}}

In this article, we analyzed the Hamiltonian pitchfork bifurcation due to changes in parameters of the potential energy. We related the qualitative changes in the topology of the potential energy surface with the changes in the linear stability of the equilibrium points. We obtained the phase space structures relevant for the transition across an index-1 saddle by computing the invariant manifolds. These are essentially the unstable periodic orbit (a normally hyperbolic invariant manifold in two DOF systems) associated with the index-1 saddle that anchors both the phase space dividing surface and the phase space conduits of transition across the index-1 saddle/s. Finally, we quantified the transition dynamics for the sub-systems (typical PES for the three cases II, III, IV) that arise due to the bifurcation by obtaining the gap time distributions and directional flux.

In future, related work will focus on correlating the statistical properties of the transition with geometry of the phase space structures. One such question involves characterizing the pulses in gap time distributions with the homoclinic and heteroclinic tangle due to the stable and unstable manifolds associated with the unstable periodic orbit. Furthermore, the quartic potential coupled with a harmonic oscillator admits a simple parametrization of the depth of the potential well which can be of use for a qualitative and quantitative analysis of the influence of depth of the PES on the reaction dynamics quantities; similar to Ref.~\cite{lyu_role_2020}.

\section*{ACKNOWLEDGMENTS}

We acknowledge the support of  EPSRC Grant No. EP/P021123/1 and Office of Naval Research (Grant No. N00014-01-1-0769).

\bibliography{bifurcation_DS,transition_dynamics}
\bibliographystyle{abbrvnat}


\end{document}